\documentclass[11pt]{article}
\usepackage[colorlinks=true, linkcolor=black]{hyperref}

\usepackage{geometry}
\usepackage{enumerate} 
\usepackage{xcolor} 
\usepackage{amsmath} 
\usepackage{nccmath} 
\usepackage{amssymb} 
\usepackage{amsthm}

\usepackage{cite,authblk}
\usepackage{titleref} 
\usepackage{caption}
\usepackage{tikz}
\usepackage{mathrsfs} 
\usepackage{longtable}
\usepackage{color}

\usetikzlibrary{quotes,angles,positioning}
\usetikzlibrary{shapes,snakes}

\usetikzlibrary{shapes.misc, positioning}

\tikzset{
	invisible/.style={opacity=0,text opacity=0},
	visible on/.style={alt=#1{}{invisible}},
	alt/.code args={<#1>#2#3}{%
		\alt<#1>{\pgfkeysalso{#2}}{\pgfkeysalso{#3}} 
	},
}

\newtheorem{theorem}{Theorem}[section]
\newtheorem{lemma}[theorem]{Lemma}
\newtheorem{corollary}[theorem]{Corollary}
\newtheorem{proposition}[theorem]{Proposition}
\newtheorem{example}[theorem]{Example}
\newtheorem{remark}[theorem]{Remark}
\newtheorem{definition}[theorem]{Definition}

\title{The Complexity of Network Satisfaction Problems for Symmetric Relation Algebras with a Flexible Atom}

\newtheorem{observation}[theorem]{Observation}

\newcommand{\bB}{\mathfrak B}

\DeclareMathOperator {\id}{Id}
\DeclareMathOperator {\Hom}{H}

\DeclareMathOperator {\Sub}{S}
\DeclareMathOperator {\Aut}{Aut}

\DeclareMathOperator {\Pol}{Pol}
\DeclareMathOperator {\NSP}{NSP}
\DeclareMathOperator {\Polc}{Pol^{can}}

\DeclareMathOperator {\age}{Age}

\DeclareMathOperator {\csp}{CSP}
\DeclareMathOperator {\nsp}{NSP}

\DeclareMathOperator {\ncp}{NCP}

\DeclareMathOperator {\rra}{RRA}

\DeclareMathOperator {\HSPf}{HSP^{fin}}
\DeclareMathOperator {\Pf}{P^{fin}}

\DeclareMathOperator {\fA}{\mathfrak{A}}
\DeclareMathOperator {\fB}{\mathfrak{B}}
\DeclareMathOperator {\fC}{\mathfrak{C}}

\author{Manuel.Bodirsky (Manuel.Bodirsky@tu-dresden.de) 
\and Simon Kn\"auer (Simon.Knaeuer@tu-dresden.de)}
\affil{Institut f\"{u}r Algebra, TU Dresden, 01062 Dresden, Germany}

\begin{document}

\maketitle

\begin{abstract}
Robin Hirsch posed in 1996 the \emph{Really Big Complexity Problem}: classify
	the computational complexity of  the network satisfaction problem for all finite relation algebras $\bf A$. 
	We provide a complete classification for the case that $\bf A$ is symmetric and has a flexible atom; {in this case,} 
	the problem is NP-complete or in P. 
	{The classification task can be reduced to the case where $\bf A$ is integral.} 
	If a finite integral relation algebra has a flexible atom, then it has a normal representation $\mathfrak{B}$. We can then study the computational complexity 
	of the network satisfaction problem of ${\bf A}$ using the universal-algebraic approach, via an analysis of the polymorphisms of $\mathfrak{B}$.
	We also use a Ramsey-type result of  Ne\v{s}et\v{r}il and Rödl and a complexity dichotomy result of Bulatov for 
	conservative finite-domain constraint satisfaction problems.  
\end{abstract}


\section{Introduction}
One of the earliest approaches to formalise  constraint satisfaction problems over infinite domains is based on \emph{relation algebras}~\cite{LadkinMaddux,HirschAlgebraicLogic}. 
We think about the elements of a relation algebra as binary relations; the algebra has operations for  intersection, union, complement, converse, and composition of relations, and constants for the empty relation, the full relation, and equality,
and is required to satisfy certain axioms. 
Important examples of relation algebras are the  Point Algebra, the Left Linear Point Algebra, 
Allen's Interval Algebra, RCC5, and RCC8, just to name a few.

{
The most important computational problems related to relation algebras that have been studied are called \emph{network satisfaction problems (NSPs)}. A finite network is given in which the  nodes are linked via the elements of a fixed relation algebra. Then the computational task of the NSP for this relation algebra is to decide whether the network satisfies a strong notion of consistency  with respect to the laws of the fixed relation algebra. 
Such NSPs for a finite relation algebra can be used
to model many computational problems in temporal and spatial reasoning  ~\cite{Duentsch,NebelRenzSurvey,Qualitative-Survey}. 
}

More than two decades ago, \cite{Hirsch} asked the \emph{Really Big Complexity Problem (RBCP)}: can we classify the computational complexity of the network satisfaction problem for every finite relation algebra? 
For example, the complexity of the network satisfaction problem for the Point Algebra and the Left Linear Point Algebra is in P~\cite{PointAlgebra,BodirskyKutzAI},
while it is NP-complete for all of the other examples mentioned above~\cite{Allen,NebelRenz}. There also exist relation algebras where the complexity of the network satisfaction problem is not in NP: Hirsch gave an example of a finite relation algebra with an undecidable network satisfaction problem~\cite{Hirsch-Undecidable}.
This result might be surprising at first sight; it
is related to the fact that the representation of a finite relation algebra by concrete binary relations over some set can be quite complicated. 
We also mention that not every finite relation algebra has a representation~\cite{LyndonRelationAlgebras}.
There are even non-representable relation algebras that are symmetric~\cite{Maddux2006}; a relation algebra is \emph{symmetric} if every element is its own converse. 


{
A simple condition that implies that a finite 
relation algebra $\mathbf{A}$ has a representation is the existence of a so-called \emph{flexible atom}~\cite{Comer1984,Maddux1982}. A flexible atom is an element of $\mathbf{A}$ that is maximally unconstrained in its interaction with other elements of the relation algebra; the formal definitions can be found in Section~\ref{section:prelim}.

 }
Such relation algebras have been studied intensively, for example in the context of the so-called \emph{flexible atoms conjecture}~\cite{Maddux1994,Alm2008}.
We will see that integral relation algebras with a flexible atom even have a \emph{normal representation},
i.e., a representation which is fully universal, square, and homogeneous~\cite{Hirsch}. 
The network satisfaction problem
for a relation algebra with a normal representation can be seen as a \emph{constraint satisfaction problem}
for an infinite structure $\mathfrak B$ 
that is well-behaved from a model-theoretic point of view; in particular, we may choose $\mathfrak B$ to be \emph{homogeneous} and \emph{finitely bounded}. 

Constraint satisfaction problems over finite domains
have been studied intensively in the past two decades, and tremendous progress has been made concerning systematic findings about their computational complexity. 
 As a highlighting result,  \cite{BulatovFVConjecture} and \cite{ZhukFVConjecture,Zhuk20}  proved  the famous Feder-Vardi dichotomy conjecture which states that every finite-domain CSP is in P or NP-complete. Both proofs build on an important connection between the computational complexity of constraint satisfaction problems and universal algebra. 

The universal-algebraic approach can also be applied to study the computational complexity of 
countably infinite 
homogeneous structures $\mathfrak B$ with finite relational signature \cite{BodirskyNesetrilJLC}. For an introduction to the field, we refer to \cite{Book}. If $\mathfrak B$ is finitely bounded,
then CSP$(\mathfrak B)$ is contained in NP (see, e.g., \cite{Bodirsky-HDR-v8}).  If $\mathfrak B$ is homogeneous and finitely bounded then a complexity dichotomy
has been conjectured, along 
with algebraic criteria that distinguish NP-complete from
polynomial-time solvable problems~\cite{BPP-projective-homomorphisms}. {The exact formulation of the conjecture from~\cite{BPP-projective-homomorphisms} in full generality requires concepts that we do not need to prove our results. In Theorem~\ref{thm:main 2} we verify these conjectures for all normal representations of finite integral symmetric relation algebras with a flexible atom, and thereby also solve Hirsch’s RBCP for symmetric relation algebras with a flexible atom.}
Phrased in the terminology of relation algebras, our result is the following. 
\begin{theorem}\label{theo:result2}
	Let $\bf A$ be a finite symmetric representable relation algebra with a flexible atom, and let $A_0$ be the set of atoms of $\bf A$. Then 
	\begin{itemize}
	\item there exists an operation 
	$f\colon A_0^6\rightarrow A_0$ 
that preserves the allowed triples of $\mathbf{A}$, satisfies 
$$\forall x_1,\ldots,x_6\in A_0.~f(x_1,\ldots,x_6)\in \{x_1,\ldots x_6\}$$
and satisfies the \emph{Siggers identity} 
	$$\forall x,y,z \in A_0.~ f(x,x,y,y,z,z)=f(y,z,x,z,x,y);$$
 in this case the network satisfaction problem for $\mathbf{A}$ is in P, or
 \item the network satisfaction problem for $\bf A$ is NP-complete.
\end{itemize}
\end{theorem}

 Moreover, {the satisfiability of the Siggers identity in Theorem~\ref{theo:result2} is} 
a decidable criterion for $\mathbf{A}$ that is a sufficient condition for the polynomial-time tractability of the network satisfaction problem of $\bf A$. We want to mention that there are several other equivalent criteria that could be used instead of the first item in the theorem, namely all characterizations of Taylor Algebras for finite conservative algebras  (see, e.g.,~\cite{Conservative}).

This article is a postprint version of an article that appeared in the Journal of Artificial Intelligence Research~\cite{BodirskyKnaeJAIR}; a conference version of the article appeared under the title ``Network satisfaction for symmetric relation algebras with a flexible
               atom'' in~\cite{BodirskyKnaeAAAI}.

{
\subsection{Proof Strategy}
Every finite integral representable relation algebra $\mathbf{A}$ with a flexible atom has a normal representation $\mathfrak{B}$; for completeness, and since we are not aware of a reference for this fact, we include a proof in Section~\ref{section:symmra with flex atom}. It follows that the classification question about the complexity of the network satisfaction problem of $\mathbf{A}$ can be translated into a question about the complexity of the constraint satisfaction problem for the relational structure $\mathfrak{B}$. 

We then associate a finite relational structure $\mathfrak{O}$  to $\mathfrak{B}$ and show that $\csp(\mathfrak{B })$ can be reduced to $\csp(\mathfrak{O})$ in polynomial-time  (Section~\ref{section:finite type structure}). 
If the structure $\mathfrak{O}$
satisfies the condition of the first case in Theorem~\ref{theo:result2}, then known results
about finite-domain $\csp$s imply that $\csp(\mathfrak{O})$ is in P~\cite{Conservative,Bulatov-Conservative-Revisited,Barto-Conservative}, 
and hence $\csp(\mathfrak{B })$ is in P, too. 
If the first case in Theorem~\ref{theo:result2}
does not apply, then 
known results about finite-domain algebras imply that there are $a,b \in A_0$ such that the \emph{canonical polymorphisms} of $\mathfrak B$ act as a projection on $\{a,b\}$~\cite{Conservative,Bulatov-Conservative-Revisited,Barto-Conservative}. 
 We first observe  NP-hardness of $\csp(\mathfrak{B })$ if $\mathfrak B$ does not have a binary injective polymorphism (Section~\ref{section: can binary injection}). If $\mathfrak B$ has a binary injective polymorphism, we use
results from structural Ramsey theory to show that $\mathfrak B$ must even have a binary injective polymorphism which is canonical (Section~\ref{sec:canon}).
This implies that none of $a,b$ equals $\id \in A$. 
We then prove that ${\mathfrak B}$ does not have a binary $\{a,b\}$-symmetric polymorphism; also in this step, we apply Ramsey theory. In Section~\ref{sec:indep} we show that this in turn implies that
all polymorphisms of ${\mathfrak B}$ must be canonical on $\{a,b\}$. 
Finally, we show that ${\mathfrak B}$ cannot have a polymorphism which acts as a majority or as a minority on $\{a,b\}$, and thus by Schaefer's theorem all polymorphisms of ${\mathfrak B}$ act as a projection on $\{a,b\}$.  This is again implied by results from Section~\ref{sec:canon}.
Finally it  follows 
 that $\csp(\mathfrak B)$ is NP-hard.  
This concludes the proof of Theorem~\ref{theo:result2}.
{
	
Our proof follows a strategy that was applied several times in the study of infinite-domain constraint satisfaction problems and recently described and generalized by \cite{MottetPinskerSmooth}. We give some details about this in Section~\ref{sec:approx}.
}

\subsection{Organisation of the Article}
Section~\ref{section:prelim} introduces all the basic concepts and tools that are used in this article. In Section~\ref{section:symmra with flex atom} we define flexible atoms and obtain first results about representable relation algebras with a flexible atom. Section~\ref{section:finite type structure} is dedicated to the atom structure and the polynomial-time tractability results. In Section~\ref{sec:ncp} we provide an additional perspective on the class of computational problems under consideration; in this section we define those problems completely without the use of the relation algebra framework. The reader can get better intuition of the class of problems studied in this article, however our results and proofs do not rely on that section. The
Sections~\ref{section: can binary injection}-\ref{sec:indep} contain the main parts of the proof as outlined in the previous paragraph.
	In Section~\ref{sec:proof} we put everything together and prove the main theorem. We end with a conclusion and a small discussion of our result.

}


\section{Preliminaries}\label{section:prelim}

We recall some basic definitions and results about universal algebra, representable relation algebras, constraint satisfaction, model theory, and structural Ramsey theory.
{
\subsection{Algebras and Structures}\label{sec: algebras and str}
{
Let $\mathbb{N}$ denote the natural numbers starting with $0$ and define $\mathbb{N}_+:= \mathbb{N}\setminus\{0\}$. 
A \emph{signature} $\tau$ is a set of relation symbols and function symbols. Each symbol is associated with a natural number, called the \emph{arity} of the symbol.  Function symbols of arity $0$ are called \emph{constant symbols}. A $\tau$-structure is a tuple $\frak A=(A; (Q^\frak A)_{Q\in \tau})$ where $A$ is a set, called the \emph{domain} of $\frak A$, such that for every $Q\in \tau$:
\begin{itemize}
	\item if $Q$ is a relation symbol of arity $n\in \mathbb{N}$, then $Q^\frak A$ is a subset of $A^n$,
	\item if $Q$ is a function symbol of arity $n\in \mathbb{N}$, then $Q^\frak A$ is an operation $A^n \rightarrow A$.
\end{itemize}}
Note that by $ A^0=\{\emptyset\}$ a subset of  $A^0$ can be seen as a Boolean value and the operation $f\colon A^0 \rightarrow A$ can be interpreted as a constant.
	{As long as there is no risk of confusion
	we will often use the function symbols for the corresponding functions, and the relation symbols for the corresponding relations, i.e., we use $Q$ instead of $Q^\frak A$.
}

Let $\frak A$  and $\frak B$ be two $\tau$-structures. A \emph{homomorphism h} from $\frak A$ to $\frak B$ is a function $h\colon A \rightarrow B$ such that
\begin{itemize}
	\item  for every relation symbol $Q$ of arity $n \in \mathbb{N}$ and every tuple $(a_1,\ldots,a_n)\in A^n$, we have that $(a_1,\ldots,a_n)\in Q^\frak A \Rightarrow  (h(a_1),\ldots,h(a_n))\in Q^\frak B$;
	\item  for every function symbol $Q$ of arity $n \in \mathbb{N}$ and every tuple $(a_1,\ldots,a_n)\in A^n$, we have that $h( Q^\frak A(a_1,\ldots,a_n) )= Q^\frak B(h(a_1),\ldots,h(a_n))$.
\end{itemize}
{In the case that $\tau$ contains only function symbols and $h$ is surjective, then $\fB$ is called \emph{homomorphic image} of $\fA$. 
In general, the homomorphism $h$  is called an \emph{embedding} if  $h$ is injective and satisfies }
\begin{itemize}
	\item for every relation symbol $Q$ of arity $n \in \mathbb{N}$ and every tuple $(a_1,\ldots,a_n)\in A^n$, we have that $(a_1,\ldots,a_n)\in Q^\frak A \Leftrightarrow  (h(a_1),\ldots,h(a_n))\in Q^\frak B$.
\end{itemize}
A surjective embedding is called an \emph{isomorphism}. An \emph{endomorphism} of  a $\tau$-structure $\fA$ is a homomorphism from $\fA$ to $\fA$ and an \emph{automorphism} of $\fA$ is an isomorphism from $\fA$ to $\fA$. We denote by $\Aut(\frak A)$ the group of all \emph{automorphisms} of $\fA$.

A $\tau$-structure $\fA$ is a \emph{substructure of  a $\tau$-structure $\fB$} if 
\begin{itemize}
	\item $A\subseteq B$;
	\item  for every relation symbol $Q$ of arity $n \in \mathbb{N}$ and every tuple $(a_1,\ldots,a_n)\in A^n$, we have that $(a_1,\ldots,a_n)\in Q^\frak A$ if and only if   $(a_1,\ldots,a_n)\in Q^\frak B$;
	\item for every function symbol $Q$ of arity $n \in \mathbb{N}$ and every tuple $(a_1,\ldots,a_n)\in A^n$, we have that $ Q^\frak A(a_1,\ldots,a_n) = Q^\frak B(a_1,\ldots,a_n)$.
\end{itemize}
Note that for every subset $A'\subseteq B$ of the domain of a $\tau$-structure $\fB$ there exists a unique substructure $\fA$ of $\fB$ with smallest domain $A$ and $A'\subseteq A$. {We call this the \emph{substructure of $\fB$ induced by $A$}.}

Let $\fA $ and $\fB$ be $\tau$ structures. The \emph{direct product} $\fC=\fA \times \fB$ is the $\tau$-structure where \begin{itemize}
	\item $A\times B$ is the domain of $\fC$;
	\item for every relation symbol $Q$ of arity $n \in \mathbb{N}$ and every tuple $((a_1, b_1),\ldots,(a_n, b_n))\in (A\times B )^n$, we have that $((a_1, b_1),\ldots,(a_n, b_n))\in Q^\frak C$ if and only if   $(a_1,\ldots,a_n)\in Q^\frak A$ and $(b_1,\ldots,b_n)\in Q^\frak B$;
	\item for every function symbol $Q$ of arity $n \in \mathbb{N}$ and every tuple $((a_1, b_1),\ldots,(a_n, b_n))\in (A\times B )^n$, we have that $ Q^\frak C ((a_1, b_1),\ldots,(a_n, b_n)) = (Q^\frak A(a_1,\ldots,a_n), Q^\frak B(b_1,\ldots,b_n) )$.
\end{itemize}
We denote the direct product $\fA\times \fA$ by $\fA^2$.  The $k$-fold product $\fA\times \cdots \times \fA$ is defined analogously and denoted by $\fA^k$.

Structures with a signature that only contains function symbols are called \emph{algebras} and structure with purely relational signature are called \emph{relational structures}. Since we do not deal with signatures of mixed type in this article, we will from now on  use the term structure for relational structures only.

}

\subsection{Representable Relation Algebras and Network Satisfaction Problems}
{
Relation algebras that are not representable have a trivial network satisfaction problem, namely the class of yes-instances of their network satisfaction problem is empty (see Definition~\ref{defi:nsp}).}
We thus omit the definition of relation algebras
and start immediately 
with the simpler definition
of \emph{representable relation algebras}; 
here we basically follow the textbook by  \cite{Maddux2006-dp}.

\begin{definition}
	Let $D$ be a set and $E \subseteq D^2$ an equivalence relation {on $D$}. \\Let  $(\mathcal{P}(E); ~\cup,  ~ \bar{}~, ~0,~1,~\id,  ~^\smile, ~\circ ) $ be an algebra 
	with the following operations:\begin{enumerate}
		\item $a \cup b := \{(x,y)  \mid (x,y)\in a \vee (x,y)\in b \} $,
		\item $\bar{a}:= E \setminus a $,
		\item $0:= \emptyset$,
		\item $1:= E$,
		\item $\id:= \{(x,x) \mid x\in D \} $,
		\item $a^\smile := \{(x,y) \mid (y,x)\in a  \}$,
		\item $a\circ b := \{(x,z) \mid \exists y \in D: (x,y)\in a  \textup{~ and~ } (y,z)\in b \}$.
	\end{enumerate}
	A subalgebra of $(\mathcal{P}(E); ~\cup,  ~ \bar{}~, ~0,~1,~\id,  ~^\smile, ~\circ ) $
	is called a \emph{proper relation algebra}.
\end{definition}

The class of \emph{representable relation algebras}, denoted by $\rra$, {consists of all algebras with signature $\tau=\{~\cup,  ~ \bar{}~, ~0,~1,~\id,  ~^\smile, ~\circ\}$ with corresponding arities $2,1,0,0,0,1$ and $2$ that are isomorphic to some proper relation algebra. }
We use bold letters (such as $\bf A$) to denote algebras from $\rra$
and the corresponding roman letter (such as $A$) to denote the domain of the algebra. 
An algebra is called 
\emph{finite} if its domain is finite.
	We call $\mathbf{A}\in\rra$ \emph{symmetric} if all its elements are symmetric, i.e., $a^\smile=a$ for every $a\in A$.

{According to the previous definition, an algebra $\mathbf{A}$ is in $\rra$ if it has an isomorphism to a proper relation algebra. Such an isomorphism is usually called the representation of $\bf A$. 
	To link the theory of relation algebras with model theory it will be convenient to view 
	representations of algebras in $\rra$ as relational structures and to not use the classical notation here. However, it is easy to see that the existence of a representation is equivalent under both definitions.
}


\begin{definition}\label{def: rep}
	Let $\mathbf{A}\in \rra$. Then a \emph{representation of $\mathbf{A}$} is a  relational structure $\mathfrak{B} $ such that 
	\begin{itemize}
		\item $\mathfrak{B}$ is an $A$-structure, i.e., the elements of $A$ are {binary} relation symbols of ${\mathfrak B}$; 
		\item The map $a \mapsto a^\mathfrak{B}$ is an isomorphism between $\mathbf{A}$ and the proper relation algebra induced by the relations of $\mathfrak B$ in
$(\mathcal{P}(1^{\mathfrak{B}}); ~\cup,  ~ \bar{}~, ~0,~1,~\id,  ~^\smile, ~\circ ) $.
		
	\end{itemize}
\end{definition}

{Let  $\mathbf{A}=(A; ~\cup,  ~ \bar{}~, ~0,~1,~\id,  ~^\smile, ~\circ ) $ be a representable relation algebra and let us define a new operation $x\cap y:= \overline{\overline{x}\cup \overline{y}}$ on the set $A$. Then the algebra  $(A; ~\cup,~\cap ,~ \bar{}~, ~0,~1 ) $ is by definition a Boolean algebra and induces therefore a partial order  $\leq$ on $A$, which is defined by
$x \leq y :\Leftrightarrow x \cup y = y$.  }
Note that for proper relation algebras, this ordering coincides with the set-inclusion order.
The minimal elements of this order in $A\setminus\{0\}$ are called \emph{atoms}. The set of atoms of $\bf A$ is denoted by $A_0$.
A tuple $(x,y,z) \in (A_0)^3$ is called an \emph{allowed triple} if $z \leq x \circ y$. 
Otherwise, $(x,y,z)$ is called a \emph{forbidden triple} { and in this case $\overline{z}\cup  \overline{x\circ y}=1 $.}
{We say that a relational $A$-structure $\mathfrak{ B}$ \emph{induces a forbidden triple (from $\mathbf{A}$)} if there exists $b_1,b_2, b_3\in B$ and $(x,y,z) \in (A_0)^3$ such that $x(b_1,b_2) , y(b_2,b_3)$ and $z(b_1,b_3)$ hold and $(x,y,z)$ is a forbidden triple. Note that a representation of $\bf A$ does not induce a forbidden triple.}
	

{
\begin{definition}\label{defi:Netw}Let $\mathbf{A}\in \rra$. An  \emph{$\mathbf{A}$-network $(V;f)$} is a finite set $V$ together with a partial function $f\colon E\subseteq V^2 \rightarrow A$, where  $E$ is the domain of $f$. 
	An $\mathbf{A}$-network $(V;f)$ is \emph{satisfiable in a representation} 
	$\mathfrak{B}$ of $\bf A$ if there exists an assignment $s\colon V \rightarrow B$ such that for all $(x,y) \in E$ the following holds: $$ (s(x), s(y))\in f(x,y)^\mathfrak{B}. $$
	An $\mathbf{A}$-network $(V;f)$ is \emph{satisfiable} if there exists a representation $\mathfrak{B}$ of $\mathbf{A}$ such that $(V;f)$ is satisfiable in $\mathfrak{B}$.
	\end{definition}


}

\begin{figure}[t]
	
				\begin{center}
				
				\begin{tabular}{|c||c|c|c|}
					\hline 
					$~\circ~$	& $~\id~$ & $~a~$ & $~b~$ \\ 
					\hline \hline
					$\id$	&  $\id$ & $a$ &$b$  \\ 
					\hline 
					$a$	&$a$  & $ \id\cup b$  & $a\cup b$ \\ 
					\hline 
					$b$	& $b$ & $a\cup b$ & $\id\cup a\cup b$ \\ 
					\hline 
				\end{tabular}

\end{center}

	\caption{Multiplication table of the representable relation algebra $\#17$.}
	\label{Tables1}
\end{figure}
{We give an example of  an instance of the $\nsp$ for a representable relation algebra.} The numbering of the representable relation algebra is by \cite{AndrekaMaddux}.
\begin{example}[An instance of $\nsp$ of representable  relation algebra $\#17$]
Let $\mathbf{A}$ be the symmetric representable relation algebra with the set of atoms $\{\id,a,b\}$ and {the values for the composition operation $\circ$ on these atoms be given by Table \ref{Tables1}. Note that  this determines the composition operation on the whole domain  of $\mathbf{A}$, which is the following set:
$$A=\{\emptyset,\id, a,b, \id \cup a,\id\cup b, a \cup b, \id\cup a \cup b\}.$$}	Let $V:=\{x_1,x_2,x_3\}$ be a set. Consider the map $f\colon V^2\rightarrow A$ given by

	\begin{align*}
	&f(x_i,x_i)=\id \text{~~for all } i\in \{1,2,3\}\\
		&f(x_1,x_2)=f(x_2,x_1)=a \\
			&f(x_1,x_3)=f(x_3,x_1)=\id\cup a \\
				&f(x_2,x_3)=f(x_3,x_2)=b\cup a. 
	\end{align*}
	The tuple $(V;f)$ is an example of an instance of $\nsp$ of $\mathbf{A} \in \rra$.
	\end{example}

We will in the following assume that for an $\mathbf{A}$-network $(V;f)$ it holds that $f(V^2)\subseteq A\setminus \{0\}$. Otherwise, $(V;f)$ is not satisfiable.
Note that every $\mathbf{A}$-network $(V;f)$  can be viewed as an $A$-structure $\mathfrak{C}$ on the domain $V$: for all $x,y \in V$ and $a\in A$ the relation $a^\mathfrak{C}(x,y)$ holds if and only if $f(x,y)=a$.

\begin{definition}\label{defi:nsp}
	The \emph{(general) network satisfaction problem} for a finite representable relation algebra $\mathbf{A}$, denoted by $\NSP(\bf A)$, is the problem of deciding whether a given  $\mathbf{A}$-network is satisfiable.
\end{definition}

\subsection{Normal Representations  and CSPs}\label{subsec: normal rep csp}

In this section we consider a subclass of $\rra$ introduced by Hirsch in 1996. For representable relation algebras $\bf A$ from this class, $\NSP(\bf A)$ corresponds naturally to a constraint satisfaction problem (CSP). In the last two decades a rich and fruitful theory emerged to analyse the computational complexity of CSPs. We use this theory to obtain results about the computational complexity of NSPs.

 In the following let $\mathbf{A} $ be in $\rra$. An $\mathbf{A}$-network $(V;f)$ is called  \emph{closed}~(transitively closed in the work by \cite{HirschAlgebraicLogic}) if  {it is total and}
	for all $x,y,z \in V$ it holds that $f(x,x)\leq \id$, $f(x,y)=f(y,x)^\smile$, and
	$f(x,z) \leq f(x,y) \circ f(y,z)$.
{It is called \emph{atomic} 
 if the range of $f$ only contains atoms from $\mathbf{A}$.}



\begin{definition}[from \cite{Hirsch}]
	Let $\mathfrak{B}$ be a representation of $\mathbf{A}$. Then $\mathfrak{B}$ is called 
	\begin{itemize}
		\item \emph{fully universal}, if every atomic closed $\mathbf{A}$-network is satisfiable in $\mathfrak{B}$;
		\item \emph{square}, if $1^\mathfrak{B}=B^2$;
		\item \emph{homogeneous}, if for every isomorphism between finite substructures of $\mathfrak{B}$ there exists an automorphism of $\mathfrak{B}$ that extends this isomorphism;
		\item \emph{normal}, if it is fully universal, square and homogeneous.
	
	\end{itemize}
	
\end{definition}

\begin{definition}\label{def: pp}
	Let $\tau$ be a relational signature. A first-order formula $\varphi(x_1,\ldots,x_n)$ is called \emph{primitive positive (pp)} if it has the  form
	$$ \exists x_{n+1},\ldots, x_{m}. (\varphi_1\wedge \cdots \wedge \varphi_s)$$
	where $\varphi_1, \ldots, \varphi_s$ are atomic formulas, i.e., formulas of the form $R(y_1,\ldots, y_l)$ for $R\in \tau$ and $y_i \in \{x_1,\ldots, x_m\}$, of the form $y=y'$ for $y,y' \in \{x_1,\ldots x_m\}$, or of the form \emph{false}. 
\end{definition}

As usual, formulas without free variables are called \emph{sentences}. {
	A relational structure is called \emph{connected} if it is not the disjoint union of two structures. 
	A \emph{connected component} of a relational structure $\mathfrak{C}$ is a substructure $\mathfrak{B}$ of $\frak C$ that is maximal with respect to domain inclusion and is connected. 
	For every primitive positive $\tau$-sentence $\varphi$  with variable set $X$  the \emph{canonical database} $D(\varphi)$ is defined as the $\tau$-structure on $X$ where $x\in X^m$ is in the relation $R^{D(\varphi)}$ for $R\in \tau$  if and only if $R(x)$ is a conjunct from $\varphi$. Conversely every relational $\tau$-structure $\mathfrak{A}$  induces a primitive positive $\tau$-sentence  on the variable set $A$, the so-called \emph{conjunctive query}, simply by taking  the conjunction all atomic formulas that hold in $\mathfrak{A}$. We say that the primitve positive $\tau$-sentences $\varphi_1,\ldots, \varphi_n$ are the \emph{connected components} of a primitive positive $\tau$-sentence $\varphi$ if they are the conjunctive queries of the connected components of the canonical database $D(\varphi)$.

	}

\begin{definition}
Let $\tau$ be a finite relational signature and 
	let $\mathfrak{B}$ be a $\tau$-structure. Then the constraint satisfaction problem of $\mathfrak{B}$ $(\csp(\mathfrak{B})) $ is the computational problem of deciding  whether a
	given primitive positive $\tau$-sentence holds in $\mathfrak B$. 
\end{definition}

{

	Consider the following translation which associates to each  $\mathbf{A}$-network $(V;f)$ a primitive positive $A$-sentences $\varphi$ as follows: the variables of  $\varphi$ are the elements of $V$ and $\varphi$ contains for every $(x,y)$ in the domain of $f$ the conjunct $a(x,y)$ if and only if $f(x,y)=a $ holds. For the other direction let $\varphi$ be an $A$-sentence with variable set $X$ and consider the $\mathbf{A}$-network $(X;f)$ with the following definition: for every $x,y\in X$, if $(x,y)$ does not appear in any conjunct from $\varphi$ we leave $f(x,y)$ undefined, otherwise let $a_1(x,y),\ldots,a_n(x,y) $ all the conjuncts from $\varphi$ that contain $(x,y)$. We compute in $\mathbf{A}$ the element $a:=a_1 \cap \ldots\cap a_n$ and define $f(x,y):=a$. 

The following theorem, which subsumes the connection between network satisfaction problems and constraint satisfaction problems  is based on this natural 1-to-1 correspondence between $\mathbf{A}$-networks and $A$-sentences.

	\begin{theorem}[\hspace{1sp}\cite{Bodirsky-HDR-v8}, see also \cite{Qualitative-Survey,BodirskyRamics}]\label{thm: nsp csp 2}
Let $\mathbf{A} \in \rra$ be finite. Then the following holds:\begin{enumerate}
	\item $\mathbf{A} $ has a representation $\fB$ such that $\nsp(\mathbf{A})$ and $\csp(\mathfrak{B})$ are the same problem up to the translation between $\mathbf{A}$-networks and $A$-sentences.
	\item If $\mathbf{A} $ has a normal representation $\fB$ the problems $\nsp(\mathbf{A})$ and $\csp(\mathfrak{B})$ are the same up to the translation between $\mathbf{A}$-networks and $A$-sentences.
\end{enumerate}
	\end{theorem}

}

%

\subsection{Model Theory}
Let $\tau$ be a finite relational signature. 
The class of finite $\tau$-structures that embed
into a $\tau$-structure $\mathfrak B$ is called the \emph{age} of $\mathfrak B$, denoted by $\age(\mathfrak{B })$. 
If $\mathcal F$ is a class of finite $\tau$-structures,
then $\textup{Forb}(\mathcal{F})$ is the class of all
finite $\tau$-structures $\mathfrak A$ such that no 
structure from $\mathcal{F}$ embeds into $\mathfrak A$. 
A class $\mathcal{C}$ of finite $\tau$ structures is called
 \emph{finitely bounded} if there exists a finite set of finite $\tau$-structures $\mathcal{F}$ 
such that ${\mathcal C} = \textup{Forb}(\mathcal{F})$ (see, e.g., \cite{MacphersonSurvey}). The structures from $\mathcal{F}$ are called \emph{bounds} or \emph{forbidden substructures}.
It is easy to see that a class $\mathcal{C}$ of $\tau$-structures is finitely bounded if and only if it is axiomatisable by a universal $\tau$-sentence.
A structure $\mathfrak{B}$ 
is called finitely bounded if 
$\age(\mathfrak B)$ is finitely bounded.



\begin{proposition}[\hspace{1sp}\cite{BodirskyRamics}]\label{prop: finite boundedness}
		Let $\mathfrak{B}$ be a normal representation of a finite  $\mathbf{A}\in \rra$. 
		Then the following holds:
		\begin{itemize}
			\item   $\mathfrak{B}$ is finitely bounded by bounds of size at most three.
	
			\item The $A_0\setminus\{\id\}$-reduct of $\mathfrak{B}$ is finitely bounded by bounds of size at most three.
		\end{itemize}

	\end{proposition}

\begin{definition}\label{def:AP}
	A class of finite $\tau$-structures has the \emph{amalgamation property} if for all structures $\mathfrak{A}, \mathfrak{B}_1, \mathfrak{B}_2 \in\mathcal{C}$ with embeddings $e_1\colon \mathfrak{A} \rightarrow \mathfrak{B}_1$ and $e_2\colon \mathfrak{A} \rightarrow \mathfrak{B}_2$ there exist a structure $\mathfrak{C}\in \mathcal{C}$ and embeddings $f_1\colon \mathfrak{B}_1\rightarrow \mathfrak{C}$ and $f_2\colon \mathfrak{B}_2\rightarrow \mathfrak{C}$ such that $f_1\circ e_1=f_2\circ e_2$.
If additionally $f_1(B_1) \cap f_2(B_2) = f_1(e_1(A)) = f_2(e_2(A))$, then we say that $\mathcal C$ has the \emph{strong amalgamation property}. 
\end{definition}

Let  $\mathfrak{B}_1, \mathfrak{B}_2 $ be $\tau$-structures. Then  $\mathfrak{B }_1\cup \mathfrak{B }_2$ is  the $\tau$-structure on the domain $B_1\cup B_2$ such that
$R^ {\mathfrak{B }_1\cup \mathfrak{B }_2} := R^\mathfrak{B }_1 \cup R^\mathfrak{B}_2$ for every $R\in \tau$. 
If Definition~\ref{def:AP} holds with 
$\mathfrak{C}:=\mathfrak{B }_1\cup \mathfrak{B}_2$ then we say that $\mathcal{C}$ has the \emph{free amalgamation property};
note that the free amalgamation property implies
the strong amalgamation property. 

\begin{theorem}[Fra\"{i}ss\'e; see, e.g., \cite{Hodges}]\label{theo:fraisse}
Let $\tau$ be a finite relational signature and let $\mathcal{C}$ be a class of finite  $\tau$-structures that is closed under taking induced substructures and isomorphisms and has the amalgamation property. Then there exists an up to isomorphism unique countable homogeneous structure $\mathfrak{B}$ such that $\mathcal{C}=\age(\mathfrak{B})$.
\end{theorem}
{

\begin{definition}\label{orbits}Let $\mathfrak{B }$ be a relational structure. A set $O\subseteq B^n$ is called an \emph{$n$-orbit} of $\Aut(\mathfrak{B })$ if $O$ is preserved by all $\alpha \in \Aut(\mathfrak{B })$ and for all $x,y \in O$ there exists $\alpha \in \Aut(\mathfrak{B })$ such that $\alpha(x)=y$.
	
\end{definition}
For a structure $\mathfrak{B }$ the automophism group $\Aut(\mathfrak{B })$ is called $\emph{transitive}$ if $\Aut(\mathfrak{B })$ has only one orbit. 
We want to remark that if $\mathfrak{B}$ is a normal representation of a finite $\mathbf{A}\in \rra$ then the 2-orbits of $\Aut(\mathfrak{B })$  are exactly the relations induced by the atoms $A_0$ of $\mathbf{A}$.

}

\subsection{The Universal-Algebraic Approach}
In this section we present basic notions for the so-called universal-algebraic approach to the study of $\csp$s.


\begin{definition}
	Let $B$ be some set. We denote by $O_B^{(n)}$ the set of all $n$-ary operations on $B$ and by $O_B := \bigcup_{n\in \mathbb{N}} O_B^{(n)}$ the set of all operations on $B$.
	A set $\mathscr{C} \subseteq O_B$ is called an \emph{operation clone} on $B$ if it contains all projections of all arities and if it is closed under composition, i.e.,  for all $f\in \mathscr{C}^{(n)} := \mathscr{C}\cap O_B^{(n)} $ and $g_1,\ldots, g_n \in \mathscr{C} \cap O_B^{(s)}$  it holds that $f(g_1,\ldots, g_n) \in \mathscr{C}$, where $f(g_1,\ldots, g_n) $  is  the $s$-ary function defined as follows
	$$f(g_1,\ldots, g_n) (x_1,\ldots, x_s ):= f( g_1(x_1,\ldots, x_s ) , \ldots, g_n(x_1,\ldots, x_s     )).$$	
\end{definition}

	An operation $f \colon B^n \to B$ is called \emph{conservative} if for all $x_1,\ldots, x_n \in B$ it holds that $f(x_1,\ldots,x_n) \in \{x_1,\ldots, x_n\}.$
A clone is called \emph{conservative} if all its operations are conservative. 
 Note that an operation clone $\mathscr{C}$ on $B$ can be considered as an algebra with domain $B$  and an infinite signature that consists of symbols for all operations in $\mathscr{C}$ (cf. Section \ref{sec: algebras and str}). In this sense a conservative  operation clone $\mathscr{C}$ on $B$ induces on every set $A\subseteq B$ a subalgebra. It is easy to see that this subalgebra corresponds to an operation clone on $A$. We call this the \emph{restriction} of $\mathscr{C}$ to $A$. 
We later need the following classical result for clones over a two-element set. 

\begin{theorem}[\hspace{1sp}\cite{Post}]
\label{post}
Let $\mathscr{C}$ be a conservative operation clone on $\{0,1\}$.
Then either  $\mathcal C$ contains only projections, or at least one of the following operations:
 \begin{enumerate}
 	\item the binary function $\min$,
 	\item the binary function $\max$,
 	\item the minority function,
 	\item the majority function.
 \end{enumerate}
\end{theorem}

Operation clones occur naturally as polymorphism clones of relational structures.
If $x_1,\dots,x_n \in B^k$ and $f \colon B^n \to B$,
then we write $f(x_1,\ldots,x_n)$ for the $k$-tuple
obtained by applying $f$ component-wise to the tuples $x_1,\ldots,x_n$.
	
\begin{definition}
Let $\mathfrak{B}$ a structure with a finite relational signature $\tau$ and let $R\in \tau$. An $n$-ary operation \emph{preserves} a relation $R^\mathfrak{B}$ if for all $x_1,\ldots, x_n \in R^\mathfrak{B}$ it holds that
	$$f(x_1,\ldots,x_n) \in R^\mathfrak{B}.$$
	If $f$ preserves all relations from $\mathfrak{B}$ then $f$ is called a \emph{polymorphism} of $\mathfrak{B}$. 
\end{definition}
In order to provide an additional view on polymorphisms we give the following definition.

\begin{definition} Let $\tau$ be a relational signature and let $\mathfrak{B }$ be a $\tau$-structure.
	For $n\in \mathbb{N}_+$ the structure $\mathfrak{B }^n$ is the $\tau$-structure on the domain $B^n$ defined as follows. Let $R\in \tau$  be of arity $l$. Then
	$$ (x^1,\ldots, x^l) \in R^{\mathfrak{B}^n} :\Leftrightarrow \forall  i\in \{1,\ldots,n\} :(x^1_i ,\ldots, x^l_i) \in R^\mathfrak{B }.$$
\end{definition}
It is easy to see that the $n$-ary polymorphisms of $\mathfrak{B}$ are precicely the homomorphisms from  $\mathfrak{B }^n$ to $\mathfrak{B}$.

The set of all polymorphisms (of all arities) of a relational structure $\mathfrak{B}$ is an operation clone on $B$, which is denoted by $\Pol(\mathfrak{B})$.
A \emph{Siggers operation} is an operation that satisfies the Siggers identity (see Theorem \ref{theo:result2}). The following result can be obtained by combining
known results from the literature.


\begin{theorem}[\hspace{1sp}\cite{Siggers,Conservative}; see also~\cite{Barto-Conservative,Bulatov-Conservative-Revisited}]\label{2elemtsubalgebrassiggers}
	Let $\bB$ be a finite structure with a finite relational signature such that $\Pol(\bB)$ is conservative. Then precisely one of the following holds:
	\begin{enumerate}
		\item There exist distinct $a,b \in B$ such that for every $f \in \Pol(\bB)^{(n)}$ the restriction
		of $f$ to $\{a,b\}^n$ is a projection. In this case, $\csp(\bB)$ is NP-complete. 
		\item $\Pol(\bB)$ contains a Siggers operation; in this case, $\csp(\bB)$ is in P. 
	\end{enumerate}
\end{theorem}

We now discuss fundamental results about
the universal-algebraic approach for constraint  satisfaction problems of structures with an infinite domain. 

\begin{theorem}[\hspace{1sp}\cite{BodirskyNesetrilJLC}]\label{polinv}
Let $\mathfrak{B}$ be a homogeneous structure with finite relational signature. Then a relation is preserved by $\Pol(\mathfrak{B})$ if and only if it is primitively positively definable in $\mathfrak{B}$.
\end{theorem}

The following definition is a preparation to formulate the next theorem which is a well-known condition that implies NP-hardness of $\csp(\mathfrak{B })$ for homogeneous structures with a finite relational signature.

\begin{definition}
	Let $\mathcal{K}$ be a class of algebras. Then we have
	\begin{itemize}
		\item $\Hom(\mathcal{K})$ is the class of homomorphic images of algebras from $\mathcal{K}$ and
		\item $\Sub(\mathcal{K})$  is the class of subalgebras of algebras from $\mathcal{K}$.
		\item 	$\Pf(\mathcal{K})$  is the class of finite products of algebras from $\mathcal{K}$.
	\end{itemize}
\end{definition}

An operation clone $\mathscr{C}$ on a set $B$ can also be seen as an algebra $\bf B$ with domain $B$ whose signature consists of the operations of $\mathscr{C}$ such that $f^{\bf B}:=f$ for all $f \in {\mathscr C}$.

The following is a classical condition for NP-Hardness, see for example Theorem 10 in the survey by \cite{BodirskySurvey}.
\begin{theorem}\label{theo:hardness} Let $\mathfrak{B}$ be a homogeneous structure with finite relational signature.
	If $\HSPf(\{\Pol(\mathfrak{B }) \}  )$ contains a 2-element algebra where all operations are projections, then $\csp(\mathfrak{B})$ is NP-hard.
\end{theorem}

In the following let $\bf A\in \rra$ be  finite and with normal representation $\mathfrak{B}$.

\begin{definition}\label{defi: 2n relation} Let $a_1,\ldots, a_n\in A_0$ be atoms of $\mathbf{A}$. Then the $2n$-ary relation $(a_1,\ldots,a_n)^\mathfrak{B}$ 
is defined as follows:
	$$(a_1,\ldots,a_n)^\mathfrak{B} := \big \{(x_1,\dots,x_n,y_1,\dots,y_n) \in B^{2n} \mid \bigwedge_{i\in \{1,\ldots,n\}} a_i^\mathfrak{B}(x_i,y_i)\big \}.$$

\end{definition}

	An operation $f\colon B^n\rightarrow B $ is called \emph{edge-conservative} if it 
	satisfies for all $x,y\in B^n$ and all $a_1,\ldots, a_n\in A_0$
	 $$(a_1,\ldots,a_n)^\mathfrak{B}(x,y)  \Rightarrow   ( f(x), f(y) )\in \bigcup_{i\in\{1,\ldots,n\}} a_i^\mathfrak{B}.$$
Note that for every $D\subseteq A_0$ the structure $\mathfrak{B }$ contains the relation $\bigcup_{a_i\in D} a_i^{\mathfrak{B}} $.
Therefore the next proposition follows immediately since polymorphisms of $\mathfrak{B}$ preserve all relations of $\mathfrak{B }$.
\begin{proposition}\label{prop:conservative}
All polymorphisms of $\mathfrak{B}$ are edge-conservative.
\end{proposition}


\begin{definition}\label{defi:canonicalfunctions}
	Let $X\subseteq A_0$.
	An operation $f\colon B^n \rightarrow B$ is called \emph{$X$-canonical (with respect to $\mathfrak{B}$)} if there exists a function $\bar{f} \colon X^n \rightarrow A_0$ such that for all $a,b\in B^n$ and $O_1, \ldots, O_n \in X$, 
if  $(a_i,b_i)\in O_i$ for all $i\in \{1,\ldots,n\}$ then $( f(a),f(b) ) \in \bar{f}(O_1, \dots,O_n)^\mathfrak{B }.$
An operation $f$ is called \emph{canonical (with respect to $\mathfrak{B}$)} if it is $A_0$-canonical. In this case we say that the behaviour $\bar{f}$ is \emph{total}. If $X\subsetneq A_0$ we call $\bar{f}$ a \emph{partial behaviour}.
	The function $\bar{f}$ is called the \emph{behaviour of $f$ on $X$}. If $X=A_0$ then
	$\bar{f}$ is just called the \emph{behaviour of $f$}.
	
\end{definition}
{We denote by $\Polc(\mathfrak{B})$ the set of all polymorphisms of $\mathfrak{B }$ that are canonical with respect to $\mathfrak{B}$.}
 It will always be clear from the context what the domain of a behaviour $\bar{f}$ is.
An operation $f \colon S^2 \to S$ is called \emph{symmetric} if for all $x,y \in S$ it holds that $f(x,y)=f(y,x) $. An $X$-canonical function $f$ is called \emph{$X$-symmetric} if the behaviour of $f$ on $X$ is symmetric.


\subsection{Ramsey Theory and Canonisation}\label{subsection:Ramsey}
We avoid giving an introduction to Ramsey theory, since the only usage of the Ramsey property is via Theorem~\ref{canonisation}, and rather refer to the survey by \cite{BodirskyRamsey} for an  introduction. 

	Let $\mathfrak{A}$ be a homogeneous $\tau$-structure such that $\age({\mathfrak A})$ has the strong amalgamation property. 
	Then the class of all $(\tau \cup \{<\})$-structures $\mathfrak A$ such that $<^{\mathfrak A}$ is a linear order and whose $\tau$-reduct (i.e. the structure on the same domain, but only with the relations that are denoted by symbols from $\tau$, see e.g.~the book by \cite{Hodges}) is from $\age({\mathfrak A})$ is a strong amalgamation class, too (see e.g.~\cite{BodirskyRamsey}). By Theorem~\ref{theo:fraisse} there exists an up to isomorphism unique countable homogeneous structure of that age, which we denote by $\mathfrak{A}_<$. 
It can be shown by a straightforward back-and-forth argument that $\mathfrak{A}_<$ is isomorphic to an expansion of $\mathfrak A$, so we identify the domain of $\mathfrak A$ and of $\mathfrak{A}_<$ along this isomorphism, and call $\mathfrak{A}_<$ the \emph{expansion of $\mathfrak A$ by a generic linear order}. 

\begin{theorem}[\hspace{1sp}\cite{NesetrilRoedlPartite,Hubicka-Nesetril-All-Those}]\label{ramseyhn}    
	Let ${\mathfrak A}$ be a relational $\tau$-structure such that $\age({\mathfrak A})$ has the free amalgamation property. Then the expansion 
	 of ${\mathfrak A}$ by a generic linear order has the Ramsey property.  
\end{theorem}
 
 The following theorem gives a connection of the Ramsey property with the existence of canonical functions and plays a key role in our analysis.

\begin{theorem}[\hspace{1sp}\cite{BP-canonical}]
\label{canonisation}
	Let $\mathfrak{B}$  be a countable homogeneous structure with finite relational signature and the Ramsey property. Let  $h\colon B^k \rightarrow B$ be an operation  and let 
	$L:=\big\{    (x_1,\dots,x_k) \mapsto \alpha ( h(  \beta_1(x_1),\ldots,\beta_k(x_k)) \mid \alpha,\beta_1\ldots,\beta_k\in \Aut(\mathfrak{B})  \big  \}.$
	
	Then there exists a canonical operation $g\colon B^k \rightarrow B$ such that for every finite $F\subset B$ there exists $g'\in L$ such that $ g' |_{F^k} =g|_{F^k} $.
\end{theorem}
\begin{remark}\label{rem:ramseyautomorp}
	Let $\mathfrak{A}$ and $\mathfrak{B}$ 
	be homogeneous structures
	with finite relational signatures. If
	$\mathfrak{A}$ and $\mathfrak{B }$ 
	have the same domain and the same automorphism group, then $\mathfrak{A}$ 
	has the Ramsey property if and only if
	$\mathfrak{B}$ has it (see, e.g., \cite{BodirskyRamsey}). 
\end{remark}


\section{Representable Relation Algebras with a Flexible Atom}
\label{section:symmra with flex atom}

%
%
%
%
%
%
%
%
%
%
%

In this section we define the concept of a flexible atom and show how to reduce the classification problem for the network satisfaction problem for a finite $\mathbf{A} \in \rra$ with a flexible atom to the situation where $\mathbf{A}$ is additionally integral (Proposition~\ref{prop:flex implies integral}).
{ A finite representable relation algebra $\mathbf{A}$ is called \emph{integral} if the element $\id$ is an atom of $\mathbf{A}$ (cf.~\cite{Maddux2006-dp}).}

 Then we show that an integral $\mathbf{A} \in \rra$ with a flexible atom has a normal representation.
 Therefore, the universal-algebraic approach is applicable; in particular, we make heavy use of polymorphisms and their connection to primitive positive definability in later sections (cf. Theorem~\ref{polinv}). Furthermore, we prove that every normal representation of a finite representable relation algebra with a flexible atom has a Ramsey expansion (Section \ref{subsec: normalre}). Therefore, the tools from Section \ref{subsection:Ramsey} can be applied, too. 
Finally we give some examples of representable relation algebras with a flexible atom (Section \ref{subsec:exam}). We start with the definition of a flexible atom.


{

\begin{definition}

Let $\mathbf{A} \in \rra$ and let $I:=\{a\in A\mid a\leq \id\}$. 
		An atom $s\in A_0\setminus I$ is called \emph{flexible} if for all $a,b \in A \setminus I$ it holds that $s\leq a\circ b$.
		
	\end{definition}
This definition can for example be found in the book by \cite[Chapter~11, Exercise~1]{HirschHodkinson}.
Note that this definition does not require the representable relation algebra $\mathbf{A}$ to be integral. This is slightly more general than the definition by \cite{Maddux1994,Maddux2006-dp}. 
As mentioned before, we show in the following section that it is sufficient for our result to classify the computational complexity of $\NSP$s for finite representable relation algebras with a flexible atom that are additionally integral.  This means  that readers who prefer this second definition by \cite{Maddux1994,Maddux2006-dp} (assuming integrality) can perfectly skip the following section and read the article with this other definition in mind. In this case representable relation algebras with a flexible atom are always implicitly integral.

}

		\subsection{Integral Representable Relation Algebras}\label{subsec: integralrel}
{ Let $\mathbf{A} \in \rra$ and let $I:=\{a\in A\mid a\leq \id\}$. The atoms  in $I\cap A_0$ are called \emph{identity atoms}. Therefore, $\mathbf{A} $ is integral if and only if  $\mathbf{A}$ has exactly one identity atom.}
{The first claim of  the following lemma is a well-known fact about atoms of representable relation algebras. This fact is used to prove a second claim about representable relation algebras with a flexible atom. In simple terms, this claim states that a representation of a non-integral representable relation algebra with a flexible atom is only ``square'' on  precisely those elements that are in one certain identity atom.}

{{
	\begin{lemma}\label{Lem:tech Rela Alger Int} Let $\mathbf{A} \in \rra$ be finite. Then  there exists for every atom $s$ a unique $ e_1\in A_0$ with $0<e_1\leq\id$  such that $s=e_1\circ s$. 
		Furthermore, if $s$ is a flexible atom then for all $e_2\in A_0$ with $0<e_2\leq\id$ and $e_2\not = e_1$ we have that $e_2\circ \overline{\id}=0$. 	
	\end{lemma}}
\begin{proof}
		Note that  $\id\circ s=s$ by definition and therefore $e\circ s\subseteq s$ for all $e\in A_0$ with $0<e\leq\id$. Since $s$ is an atom either $e\circ s=0$  or $e\circ s=s$. { By $\id=\bigcup \{e\in A_0\mid 0<e\leq\id\}$ and $\id \circ s=s$ there exists at least one $0<e\leq\id$ with  $e \circ s=s$.  In the next step we prove uniqueness of such an element $e$.
		 Assume for contradiction that there exist distinct $e_1,e_2\in A_0$ with $0<e_1\leq\id$ and $0<e_2\leq\id$ such that $e_1\circ s=s$ and $e_2\circ s=s$. Note that $e_1\circ e_2=0$ since $e_1$ and $e_2$ are identity atoms. Therefore, we get $$0=0\circ s= (e_1\circ e_2)\circ s= e_1\circ( e_2\circ s)= e_1\circ s=s,$$which is a contradiction since $s $ is an atom.}
	  This proves the first statement.
		
		For the second statement assume  for contradiction that there exists  $e_2\in A_0\setminus\{e_1\}$ such that $e_2\leq\id$ and $e_2\circ \overline{\id} \not = 0$. Let $a$ be an atom with $a\leq e_2\circ\overline{\id} $. {Since $e_1\circ e_2=0$  we get $e_1\circ a\leq e_1\circ( e_2\circ\overline{\id}) =    (e_1\circ e_2)\circ\overline{\id} = 0 \circ \overline{\id} =0$. Since $s$ is a flexible atom it holds that 
			$s\leq a\circ a^\smile$ and therefore 
			$$s=e_1\circ s\leq e_1 \circ (a \circ a^\smile )   = (e_1 \circ a) \circ a^\smile    =0\circ a^\smile =0,$$ which is a contradiction. }
	 \end{proof}

}

{
\begin{proposition}\label{prop:flex implies integral}
	Let $\mathbf{A} \in \rra$ be finite and with a flexible atom. Then there exists a finite integral $\mathbf{A}' \in \rra$ with a flexible atom such that the following statements hold:
	\begin{enumerate}
	{\item 	There exists a polynomial-time many-one  reduction from $\nsp(\mathbf{A})$ to $\nsp(\mathbf{A}')$.}
		\item There exists a polynomial-time many-one reduction from $\nsp(\mathbf{A}')$ to $\nsp(\mathbf{A})$.
		\item The atom structure of  $\mathbf{A} $ has a polymorphism that satisfies the Siggers identity if and only if  the atom structure of $\mathbf{A}' $ has such a polymorphism (see Definition~\ref{def: atom str}).
	\end{enumerate}
\end{proposition}
\begin{proof}If $\bf A$ is integral there is nothing to be shown. So assume that $\mathbf{A}$ is not integral and let $s$ be a flexible atom. Let $\mathfrak{B }$ be a representation of $\mathbf{A}$ such that $\nsp(\mathbf{A})$ and $\csp(\mathfrak{B})$ are the same problem up to the translation between $\mathbf{A}$-networks and $A$-sentences. Such a representation exists by Theorem~\ref{thm: nsp csp 2}. Let $(x,y)\in \overline{\id}^\mathfrak{B }$ and let $e_1\in A_0$ be the unique element with $e_1\leq\id$ and $s=e_1\circ s$  that exists by  Lemma \ref{Lem:tech Rela Alger Int}. The second statement of  Lemma \ref{Lem:tech Rela Alger Int} implies $e_1\circ \overline{\id}=\overline{\id}$ and therefore we have that $(x,x)\in e_1^\mathfrak{B }$ and $(y,y)\in e_1^\mathfrak{B }$.	Let $\mathfrak{C}'$ be the substructure of $\mathfrak{B }$ on the domain $\{x\in B\mid (x,x)\in e_1^\mathfrak{B }\}$. The set of relations of $\mathfrak{C }'$ clearly induces a proper relation algebra which is integral. We denote this representable relation algebra by $\mathbf{A}'$.  Note that we can also consider $A'$ as a subset of $A$.
	Let $\mathfrak{B' }$ be the representation of $\mathbf{A}'$ such that $\nsp(\mathbf{A}')$ and $\csp(\mathfrak{B}')$ are the same problem up to the translation between $\mathbf{A}'$-networks and $A'$-sentences. As before, such a representation exists by Theorem~\ref{thm: nsp csp 2}.
	
	{Proof of 1.: Note that if  a connected instance of $\csp(\fB)$ is  satisfiable, then either all variables are mapped to an atom from the subset of $A_0$ that corresponds to $A'_0$ or  all variables are mapped to one element $x$ with $e^*(x,x)$ and $e^*\in \{e\in A_0 \mid e\leq \id \text{~and~}e\not = e_1\}$. This leads to the following polynomial-time many-one reduction from $\csp(\fB)$ to $\csp(\fB')$, which proves together with  Theorem~\ref{thm: nsp csp 2} the claim that there exists a polynomial-time many-one reduction from $\nsp(\mathbf{A})$ to $\nsp(\mathbf{A}')$.
		Consider the following algorithm: For a given primitive positive $A$-sentence $\varphi$ it computes
		the  connected components $\varphi_1,\ldots,\varphi_n$ (see paragraph after Definition~\ref{def: pp}). Then it checks  for every $i \in \{1,\ldots,n\}$ whether there exists $e^*\in \{e\in A_0 \mid e\leq \id \text{~and~}e\not = e_1\}$ such that for every conjunct $a(x,y)$ of $\varphi_i$  it holds that $e^*\leq a$. Let $I\subseteq \{1,\ldots, n\}$ be the set of indices for which this is not the case. 
			Then the algorithm  defines new $A'$-sentences $\varphi_i'$ from $\varphi_i $ for every $i\in I$ by replacing  every conjunct $a(x,y) $ with $a'(x,y)$ where $a':=a\setminus \{e\in A_0 \mid e\leq \id \text{~and~}e\not = e_1\}$. 
			 Let $\varphi'$ be the conjunction of all $\varphi_i'$ for $i\in I$. 
			
			It is easy to see that $\varphi'$ is  computable from $\varphi$ in polynomial time and that  $\varphi$ is a satisfiable  instance of $\csp(\fB)$  if and only if $\varphi'$ is a satisfiable instance of $\csp(\fB')$.
	}

	Proof of 2.: 
	Consider now an $\mathbf{A}'$-network $(V;f')$.
	We claim that $(V;f')$  is satisfiable as an $\mathbf{A}'$-network if and only if it is satisfiable as an $\mathbf{A}$-network.
		Suppose that $(V;f') $ is satisfiable in a representation $\mathfrak{D}'$ of $\mathbf{A}'$ by an assignment $\alpha$.  Let $y_i$ be fresh elements for every atom $e_i\leq \id$ with $e_i\not = e_1$.  We build the disjoint union of $\mathfrak{D }'$ with one-element $\{e_i\}$-structures $(\{y_i\}; \{(y_i,y_i)\} ) $  and then close the structure under union and intersection of binary relations. This results in a representation of $\mathbf{A}$ that satisfies $(V;f') $ again by the assignment $\alpha$. For the other direction, if $(V;f') $ is satisfiable in a representation $\mathfrak{D}$ of $\mathbf{A}$ we can again consider the substructure on the domain $(x,x)\in e_1^\mathfrak{D }$ and get a representation of $\mathbf{A}'$ that satisfies $(V;f')$.

	Proof of 3.: Let $g$ be a polymorphism of the atom structure of  $\mathbf{A} $  that satisfies the Siggers identity.  By assumption $g$  satisfies $$\forall x_1,\ldots,x_6\in A_0.~g(x_1,\ldots,x_6)\in \{x_1,\ldots x_6\}$$ and therefore the restriction of $g$ to $(A_0\cap A')^6$  is a polymorphism of the atom structure of  $\mathbf{A}' $ .
	
	For the other direction choose an arbitrary ordering of the atoms $\{e\in A_0 \mid e\leq \id \text{~and~}e\not = e_1\}= \{l_1,\ldots,l_j\}$. If $g$ is a Siggers polymorphism of the atom structure of  $\mathbf{A}' $  one can extend $g$ to an operation $g^*\colon A^6\rightarrow A$ by defining
	$$
	g^*(x_1,\ldots,x_6):= \left\{
	\begin{array}{ll}
		\min(\{l_1,\ldots,l_j \} \cap \{x_1,\ldots, x_6\})& \text{~if~}\{l_1,\ldots,l_j \} \cap \{x_1,\ldots, x_6\}\not = \emptyset,\\
		g (x_1,\ldots,x_6) & \, \text{otherwise.} \\
	\end{array}
	\right. 
	$$
	It is easy to see that this operation satisfies also the Siggers identity. Furthermore, since every atom $e$ from $\{e\in A_0 \mid e\leq \id \text{~and~}e\not = e_1\}$ is only contained in allowed triples of the form $(e,e,e)$ it follows that $f^*$ preserves the allowed triples from $\mathbf{A}$ (see after Definition~\ref{def: atom str}).

\end{proof}
}

\subsection{Normal Representations}\label{subsec: normalre}

	Let  $\mathbf{A \in \rra}$ be for the rest of the section finite, integral, and with a flexible atom $s$. 
We consider the following subset of $A$:
$$A-s:= \{  a\in A \mid s\not \leq a \}.$$
{Let $(V,g)$ be an $\mathbf{A}$-network and let $\mathfrak{C}$ be the corresponding $A$-structure (see paragraph before Definition~\ref{defi:nsp}). }
Let $\mathfrak{C}-s$ be the $(A-s)$-structure on the same domain $V$ as $\mathfrak{C}$ such that for all $x,y \in V$ and $a\in (A - s) \setminus \{0\}$  we have
$$ a^{\mathfrak{C}-s}(x,y)  \text{~  if and only if ~} ( a^{\mathfrak{C}}(x,y) \vee (a\cup s)^{\mathfrak{C}}(x,y)).$$
We call $\mathfrak{C}-s$ the \emph{$s$-free companion} of an $\mathbf{A}$-network $(V,f)$.

The next lemma  follows directly from the definitions of flexible atoms and $s$-free companions.

\begin{lemma}\label{lemma:free amalgamation}
Let $\mathcal{C}$ be the class of $s$-free companions of atomic closed $\mathbf{A}$-networks.
	Then $\mathcal{C}$ has the free amalgamation property.
\end{lemma}
{\begin{proof}Let  $\mathfrak{A}, \mathfrak{B}_1,$ and $\mathfrak{B}_2$ be structures in  $\mathcal{C}$ with embeddings $e_1\colon \mathfrak{A} \rightarrow \mathfrak{B}_1$ and $e_2\colon \mathfrak{A} \rightarrow \mathfrak{B}_2$. Since $s$ is a flexible atom and $\mathcal{C}$ is a class of 
		$s$-free companions we get that the structure	$\mathfrak{C}:=\mathfrak{B }_1\cup \mathfrak{B}_2$ is in $\mathcal{C}$. Therefore, the natural embeddings $f_1\colon \mathfrak{B}_1\rightarrow \mathfrak{C}$ and $f_2\colon \mathfrak{B}_2\rightarrow \mathfrak{C}$ prove the free amalgamation property of $\mathcal{C}$. \end{proof}}

As a consequence of this lemma we obtain the following.
\begin{proposition}\label{prop: flex atom implies normal rep}
 $\mathbf{A}$ has a normal representation $\mathfrak{B}$. 
\end{proposition}
\begin{proof}
	Let $\mathcal{C}$ be the class from Lemma \ref{lemma:free amalgamation}. This class is closed under taking substructures and isomorphisms. By Lemma \ref{lemma:free amalgamation} it also has the amalgamation property and therefore we get by Theorem~\ref{theo:fraisse} a homogeneous structure $\mathfrak{B}'$ with $\age(\mathfrak{B}') =\mathcal{C}$. Let $\mathfrak{B}''$ be the expansion of $\mathfrak{B}'$ by the following relation
	$$  s(x,y) :\Leftrightarrow \bigwedge_{ a\in A_0\setminus \{s\} }  \neg a^{\mathfrak{B}'}(x,y).$$
	Let $\mathfrak{B}$ be the (homogeneous) expansion of $\mathfrak{B}''$ by all Boolean combinations of relations from ${\mathfrak B}''$. 
	Then ${\mathfrak B}$ is a representation of  $\mathbf{A}\in \rra$. 
	Since $\age(\mathfrak{B}')$ is the class of all atomic closed $\mathbf{A}$-networks, 
	$\mathfrak{B}$ is fully universal. 
	The definition of $s$ witnesses  that  
	$\mathfrak{B}$ is a square representation of 
	${\bf A}$: for all elements $x,y\in B$ there exists an atom $a\in A_0$ such that $a^\mathfrak{B}(x,y)$  holds.
	 \end{proof}

The next theorem is another consequence of  Lemma \ref{lemma:free amalgamation}.

\begin{theorem}\label{theorem:ramsey order}Let $\mathfrak{B}$ be a normal representation of  $\mathbf{A}$. Let $\mathfrak{B}_<$ be the expansion of $\mathfrak{B}$ by a generic linear order. Then $\mathfrak{B}_<$ has the Ramsey property. 
\end{theorem}
\begin{proof}
	Let $\mathfrak{B }'$ be the $(A_0 \setminus\{s\})$-reduct of $\mathfrak{B}$. The age of this structure has the free amalgamation property by Lemma~\ref{lemma:free amalgamation}. Therefore,  Theorem~\ref{ramseyhn} implies that the expansion of $ \mathfrak{B }'$ by a generic linear order has the Ramsey property. By Remark~\ref{rem:ramseyautomorp} the structure $\mathfrak{B}_<$  also has the Ramsey property since $\mathfrak{B}_<$ and $(\mathfrak{B}')_<$ have the same automorphism group.	 \end{proof}

\begin{remark}\label{rem:proper RA Bor} The binary first-order definable relations of $\mathfrak{B}_<$ form a proper relation algebra since $\mathfrak{B}_<$ has quantifier-elimination (see \cite{Hodges}). By the definition of the generic order the atoms of this proper relation algebra are of the following form
		\begin{itemize}
			\item $a^{\mathfrak{B}_<} \cap <^{\mathfrak{B}_<}$ for $ a\in A_0\setminus\{\id\} $, or 
			\item $a^{\mathfrak{B}_<} \cap >^{\mathfrak{B}_<} $ for $ a\in A_0\setminus\{\id\} $, or
			\item $\id$,
		\end{itemize}
	where the relation $>^{\mathfrak{B}_<} $ consists of all tuples $(x,y)$ such that $(y,x)\in <^{\mathfrak{B}_<}$ holds.
	\end{remark}
%

\subsection{Examples}\label{subsec:exam}

\begin{figure}[t]
	\begin{center}
		
		\begin{minipage}[b]{63mm}
			\begin{center}
				
				\begin{tabular}{|c||c|c|c|}
					\hline 
					$~\circ~$	& $~\id~$ & $~a~$ & $~b~$ \\ 
					\hline \hline
					$\id$	&  $\id$ & $a$ &$b$  \\ 
					\hline 
					$a$	&$a$  & $\{\id,a,b\}$  & $\{a,b\}$ \\ 
					\hline 
					$b$	& $b$ & $\{a,b\}$ & $\{\id,a,b\}$ \\ 
					\hline 
				\end{tabular} 
				
			\end{center}
		\end{minipage}
		\begin{minipage}[b]{63mm}
			\begin{center}
				
				\begin{tabular}{|c||c|c|c|}
					\hline 
					$~\circ~$	& $~\id~$ & $~a~$ & $~b~$ \\ 
					\hline \hline
					$\id$	&  $\id$ & $a$ &$b$  \\ 
					\hline 
					$a$	&$a$  & $ \{\id,b\}$  & $\{a,b\}$ \\ 
					\hline 
					$b$	& $b$ & $\{a,b\}$ & $\{\id,a,b\}$ \\ 
					\hline 
				\end{tabular} 
			\end{center}
		\end{minipage}
	\end{center}
	\caption{Multiplication tables of representable relation algebras $\#18$ (left) and $\#17$ (right).}
	\label{Tables}
\end{figure}
We give two concrete examples of finite integral symmetric representable relation
algebras with a  flexible atom (Examples \ref{exam:rado} and \ref{exam:henson}), and a systematic way of building such algebras from arbitrary representable relation algebras (Example \ref{examp:+fexible atom}).
The numbering of the algebras in the examples is by  \cite{AndrekaMaddux}.

\begin{example}[Representable relation algebra $\#18$] \label{exam:rado}
	The representable relation algebra $\#18$ has three atoms, namely the identity atom $\id$ and two symmetric atoms $a$ and $b$. The multiplication table for the atoms is given in Fig.~\ref{Tables}. In this representable relation algebra the atoms $a$ and $b$ are flexible. 
	Consider  the countable, homogeneous, undirected graph  $\mathfrak{R}=(V; E^\mathfrak{R})$, whose age is the class of all finite undirected graphs (see, e.g., \cite{Hodges}), also called  the \emph{Random graph}. The expansion of $\mathfrak{R}$ by all binary first-order definable relations is a normal representation of the algebra $\#18$. In this representation the atoms $a$ and $b$ are interpreted as the relation $E^\mathfrak{R}$ and the relation $N^\mathfrak{R}$, where $N^\mathfrak{R}$ is defined as $\neg E(x,y) \wedge x\not=y$.
	
\end{example}

\begin{example}[Representable relation algebra  $\#17$]\label{exam:henson}
The representable relation algebra $\#17$ also consists of three symmetric atoms. The multiplication table in Fig.~\ref{Tables} shows that in this algebra the element $b$ is a flexible atom. To see that $a$ is not a flexible atom, note that $a\not\leq a\circ a= \{\id,b\} $.
	Let  $\mathfrak{N}=(V; E^\mathfrak{N})$ be the countable, homogeneous, undirected graph, whose age is the class of all finite undirected graphs that do not embed the complete graph on three vertices (see, e.g., \cite{Hodges}). This structure is  called a \emph{Henson graph}.
		If we expand $\mathfrak{N}$ by all binary first-order definable relations we get a normal representation of the algebra $\#17$. To see this note that we interpret $a$ as the relation $E^\mathfrak{N}$. That $\mathfrak{N}$ is triangle free, i.e. triangles of $E^\mathfrak{N}$ are forbidden, matches with the fact that $a\not\leq a\circ a$ holds in the representable relation algebra.
	
\end{example}

\begin{example}\label{examp:+fexible atom}
	
	Consider an arbitrary	finite, integral $\mathbf{A}=(A; \cup, \bar{}, 0,1, \id,^\smile,\circ) \in \rra$. Clearly $\mathbf{A}$ does not have a flexible atom $s$ in general. 
	Nevertheless  we can expand the domain of $\mathbf{A}$ to implement an ``artificial'' flexible atom.
	
	Let $s$ be some symbol not contained in $A$. 
	Let us mention that every element in $\mathbf{A}$ can uniquely be written as a union of atoms from $A_0$. 
	Let $A'$ be the set of all subsets of $A_0\cup \{s\}$. The set $A'$ is the domain of our new algebra $\mathbf{A}' $. Note that on $A'$ there exists the subset-ordering  and $A'$ is closed under set-union and complement (in $A_0\cup \{s\}$)
	We define $s$ to be symmetric and therefore get the following unary function $^*$  in $\mathbf{A}'$ as follows. For an element $x\in A'$ we define
	$$x^* := \left\{
	\begin{array}{ll}
	y^\smile \cup \{s\} & \text{~if~}x = y\cup \{s\}  \text{~for~} y\in A,\\
	x^\smile& \, \text{otherwise.} \\
	\end{array}
	\right. $$
	The new function symbol $\circ_A' $  in  $\mathbf{A}'$ is defined on the atoms $A_0\cup \{s\}$ as follows:
	$$x\circ_{A'}y := \left\{
	\begin{array}{ll}
	A_0\cup \{s\}& \text{~if~} \{s\}= \{x,y\},\\
	(A_0\setminus \{\id\}) \cup \{s\}& \text{~if~}\{s,a\}= \{x,y\} \text{~for~} a\in A_0\setminus\{s,\id\},\\
	\{a\}& \text{~if~}\{\id,a\}= \{x,y\}\text{~for~} a\in A_0\cup\{s\},\\
	(x\circ y) \cup \{s\} & \, \text{otherwise.} \\
	\end{array}
	\right. $$
	One can check that $\mathbf{A}'=(A'; \cup, \bar{}, \emptyset,A_0\cup\{s\}, \id,^*,\circ_{A'})$ is a finite integral representable relation algebra with a flexible atom $s$. Note that the forbidden triples of $\mathbf{A}'$ are exactly those of $\mathbf{A}$ together with triples which are permutations of $(s,a,\id)$ for some $a\in A_0$.

\end{example}

\section{Polynomial-time Tractability}\label{section:finite type structure}
In this section we introduce for every finite 
$\mathbf{A}\in \rra$ 
an associated finite structure, called the \emph{atom structure} of $\bf A$. Note that it is closely related, but not the same, as the type structure introduced by \cite{Bodirsky-Mottet}. In the context of relation algebras the atom structure has the advantage that its domain is the set of atoms of $\mathbf{A}$, rather than the set of 3-types, which would be the domain of the type structure of~\cite{Bodirsky-Mottet}; 
hence, our domain is smaller and has some advantages  on which the main result of this section (Proposition~\ref{tractabilityB}) is based.
{Up to a some differences in the signature, our atom structure is the same as the atom structure introduced by~\cite{LyndonRelationAlgebras} which was used there for different purposes (see also~\cite{Maddux1982,HirschHodkinsonStrongly,HirschJacksonK}).}

Let $\mathfrak{B}$ be a normal representation of a finite $\mathbf{A}\in \rra$. We will reduce $\csp(\mathfrak{B})$ to the $\csp$ of the atom structure of $\bf A$. This means that  if the $\csp$ of the atom structure is in P, then
so are $\csp(\mathfrak{B})$ and $\nsp(\bf A)$. For our main result we will show later that every network satisfaction problem for a finite integral symmetric representable relation algebra with a flexible atom that cannot be solved in polynomial time by this method 
is NP-complete. 



\begin{definition}\label{def: atom str}
The \emph{atom structure of $\bf A \in \rra$} is the finite relational structure $\mathfrak{O}$ with domain $A_0$ and the following relations:

	\begin{itemize}
		\item for every $x \in A$ the unary relation $x^\mathfrak{O}:= \{a\in A_0 \mid   a\leq x \}$,
		\item the binary relation  $E^\mathfrak{O}:= \{ (a_1,a_2)   \in A_0^2  \mid  a_1^\smile=a_2 \}$,
		\item the ternary relation $H^\mathfrak{O}:=\{ (a_1,a_2,a_3) \in A_0^3\mid  a_3\leq a_1\circ a_2  \}$.
	\end{itemize}
\end{definition}

Note that the relation $H^\mathfrak{O}$ consists of the allowed triples of  $\bf A \in \rra$. We say that an operation \emph{preserves the allowed triples} if it preserves the relation $H^\mathfrak{O}$.

\begin{proposition}  \label{reductionBtoO}
		Let $\mathfrak{B}$ be a fully universal representation of a finite $\mathbf{A}\in \rra$.
	There is  a polynomial-time reduction from 
	$\csp(\mathfrak{B})$ 
	to $\csp(\mathfrak{O})$.
\end{proposition}

\begin{proof}
	Let $\Psi$ be an instance of $\csp(\mathfrak{B} )$ with variable set $X=\{x_1,\ldots,x_n\}$. 
	We construct an instance $\Phi$ of $\csp(\mathfrak{O})$ as follows.
	The variable set $Y$ of $\Phi$ is given by
	$Y:= \{(x_i,x_j)\in X^2 \mid i \leq j  \}$.
	The constraints of $\Phi$ are of   the two kinds:
	\begin{enumerate}
		\item Let $a\in A$ be an element of $\mathbf{A}\in \rra$ and let $a((x_i,y_j))$ be an atomic formula of $\Psi$. If $i< j$, then we add the atomic (unary) formula $a((x_i,x_j) )$ to $\Phi$; otherwise we add the atomic formula $a^\smile((x_j,x_i))$. If $j=i$ we additionally add $\id((x_i,x_j))$. 
		
		\item Let $x_i,x_j, x_l\in X$ be such that $i\leq j\leq l$. Then we add the atomic formula \\$H( (x_i,x_j),(x_j,x_l),(x_i,x_l) ) $ to $\Phi$.
	\end{enumerate}
		It remains to show that this reduction is correct. Let $s\colon X \rightarrow B$ be a satisfying assignment for $\Psi$. This assignment maps every pair of variables $x_i$ and $x_j$ to a unique atom in $A_0$ and therefore induces a map $s' \colon Y \rightarrow A_0$. The map $s'$  clearly satisfies all atomic formulas introduced by (1.). To see that it also satisfies all formulas introduced by (2.) note that $s$ maps the elements $x_i,x_j, x_l\in X$ to a substructure of $\mathfrak{B}$, which does not induces a forbidden triple.
		
		For the other direction assume that $s'\colon Y \rightarrow A_0$ is a satisfying assignment for $\Phi$. This induces an $A$-structure $\mathfrak{X}$ on $X$ (maybe with some identification of  variables)  as follows: we add $(x_i,x_j)$ to the relation $a^\mathfrak{X}$ if  $i\leq j$ and  $s'((x_i,x_j))=a$; if otherwise $j<i$ and $s'((x_j,x_i))=a$ we add $(x_i,x_j)$ to the relation $(a^\smile)^{\mathfrak{X}}$.
		 It is clear that no forbidden triple from $\mathbf{A}$ is induced by $\mathfrak{X}$. Also note that $\mathfrak{X}$ satisfies $\Psi$ by the choice of the (unary) constraints of the first kind. Since $\mathfrak{B}$ is a fully universal 
		representation the structure $\mathfrak{X}$ is a substructure of $\mathfrak{B}$. Hence, the instance $\Psi$ is satisfiable in $\mathfrak{B}$. 	
		 \end{proof}


The atom structure has another property which is fundamental for
our proof of Theorem 1.
Recall that every canonical polymorphism $f$ induces a behaviour $\bar{f} \colon A_0^n \rightarrow A_0$. In the next proposition we show that then $\bar{f}$ is a polymorphism of $\mathfrak{O}$. Moreover the other direction also holds. Every $g\in \Pol(\mathfrak{O})$ is the behaviour of a canonical polymorphism of $\mathfrak{B}$.

\begin{proposition}\label{prop:canon poly=type struc}
	Let $\mathfrak{B}$ be a normal representation of a finite $\mathbf{A}\in \rra$. 
	\begin{enumerate}
		\item Let $g\in \Pol(\mathfrak{B})^{(n)}$ be canonical and let $\overline{g} \colon A_0^n\rightarrow A_0$ be its behaviour. Then $\overline{g} \in \Pol(\mathfrak{O})^{(n)}$. 
		\item Let $f\in \Pol(\mathfrak{O})^{(n)}$. Then there exists a canonical $g\in \Pol(\mathfrak{B})^{(n)}$ whose behaviour equals $f$. 
	\end{enumerate}
\end{proposition}

\begin{proof}
	For (1): Let $g\in \Pol(\mathfrak{B})^{(n)}$ be canonical and let $c^1,\ldots,c^n \in H^\mathfrak{O}$. Then by the definition of $H^\mathfrak{O}$ there exist  tuples $x^1,\ldots,x^n \in B^3$ such that for all $i\in\{1,\ldots,n\}$ we have $${c_1^i}^\mathfrak{B} (x^i_1,x^i_2), ~ {c_2^i}^\mathfrak{B} (x^i_2,x^i_3), \textup{~and~} {c_3^i}^\mathfrak{B} (x^i_1,x^i_3).$$
	
	We apply the canonical polymorphism $g$ and get $y:=g( x^1,\ldots,x^n )\in B^3$. Then there exists an allowed triple $(d_1,d_2,d_3)\in A_0^3$ such that
	$$ d_1^\mathfrak{B }(y_1,y_2), ~ d_2^\mathfrak{B }(y_2,y_3),  \text{~and~} d_3^\mathfrak{B }(y_1,y_3).$$
	
	We have that $d=(d_1,d_2,d_3) \in H^\mathfrak{O}$ and by the definition of the behaviour of a canonical function we get $\overline{g}(c^1,\ldots,c^n)=d$. 
	The other relations in $\mathfrak{O}$ are preserved trivially and therefore $\overline{g} \in \Pol(\mathfrak{O})^{(n)}$ .
	
	For (2): Since  $\mathfrak{B }$ is fully universal and homogeneous it follows by a compactness argument (see e.g.,  Lemma 2 by~\cite{BodDalJournal}) that every
	countable $A_0$-structure which does not induce a forbidden triple and is square has a homomorphism to $\mathfrak{B }$. It is therefore enough to show that every operation $h\colon B^n\rightarrow B$ with behaviour $f$   does not induce a forbidden triple in the image. 
	Let  $x^1,\ldots,x^n \in B^3$ be such that the application of a canonical function with behaviour $f$ on $x^1,\ldots,x^n $ would give a tuple $y\in B^3$ with $d=(d_1,d_2,d_3)\in A_0^3$ such that
	$$ d_1^\mathfrak{B }(y_1,y_2), \; d_2^\mathfrak{B }(y_2,y_3),  \text{~and~} d_3^\mathfrak{B }(y_1,y_3).$$
	Since $f$ preserves $H^\mathfrak{O}$ the triple $d$ is not forbidden.
	 \end{proof}

Recall from Proposition~\ref{prop:conservative} that polymorphisms of $\mathfrak{B}$ are edge-conservative. Note that this implies that polymorphisms of $\mathfrak{O}$ are conservative.  In fact, Theorem \ref{2elemtsubalgebrassiggers} and the Proposition~\ref{prop:canon poly=type struc} imply the following.

\begin{proposition}\label{tractabilityB}
	If $\Pol(\mathfrak{B})$ contains a canonical polymorphism $s$ whose behaviour $\overline{s}$ is a Siggers operation in $\Pol(\mathfrak{O}) $ then $\csp(\mathfrak{B})$ is in P.
\end{proposition}



We demonstrate how this result can be used to prove polynomial-time tractability of  $\nsp(\mathbf{A})$ for a symmetric, integral $\mathbf{A}\in \rra$ with a flexible atom.

\begin{example}[Polynomial-time tractability the NSP of representable relation algebra  $\#18$]\quad
	The polynomial-time tractability of the NSP of the representable relation algebra  $\#18$ (see Example~\ref{exam:rado})
		was first shown by  \cite{HirschCristiani} (see also Section 8.4 of \cite{BodPin-Schaefer-both}).
	Here we consider the following function $\bar{s}\colon \{\id,a,b\}^6\rightarrow \{\id,a,b\}$.
	
	$$\bar{s}(x_1,\ldots,x_6) := \left\{
	\begin{array}{ll}
	a& \text{~if~} a\in \{x_1,\ldots,x_6\},\\
	b& \text{~if~} b\in \{x_1,\ldots,x_6\} \text{~and~} a\not \in \{x_1,\ldots,x_6\},\\
	\id& \text{~otherwise.} 
	\end{array}
	\right. $$
	Let $\mathfrak{R}'$ be the normal representation of the algebra $\#18$ given in Example \ref{exam:rado}. 
	Note that $\bar{s}$ is the behaviour of an injective, canonical polymorphism of $\mathfrak{R}$. The injectivity  
	follows from the last line of the definition; if $\bar{s}(x_1,\ldots,x_6)=\id$ then $\{x_1,\ldots,x_6\}=\{\id\}$. 
	Therefore $\bar{s}$ preserves all allowed triples, since in the algebra $\#18$ the only forbidden triples involve $\id$.
	One can check that $\bar{s}$ is a Siggers operation and therefore we get by Proposition \ref{tractabilityB}  that $\nsp(\#18)$ is in P.
	\end{example}

\begin{example} Consider the construction of representable relation algebras with a flexible atom from Example \ref{examp:+fexible atom}. It is easy to see that  $\nsp(\mathbf{A})$ for a finite integral $\mathbf{A}\in\rra$ has a polynomial-time reduction to $\nsp(\mathbf{A'} )$ where $\mathbf{A'}$ is the representable relation algebra with a flexible atom that is constructed in Example \ref{examp:+fexible atom}. We get as a consequence that if a normal representation of $\mathbf{A'}$ satisfies the condition of Proposition \ref{tractabilityB} then $\nsp(\mathbf{A})$ is in P.
	\end{example}

 Theorem \ref{2elemtsubalgebrassiggers} has  additionally to Proposition \ref{tractabilityB} another important consequence.
\begin{corollary}\label{trivialsubalgebraO}
	If $\Pol(\mathfrak{O})$ does not have a Siggers operation then there exist elements $a_1,a_2\in A_0$ such that the restriction of every operation from $\Pol(\mathfrak{O})^{(n)}$ to $\{a_1,a_2\}^n$ is a projection. 
\end{corollary}




\section{Network Consistency Problems }\label{sec:ncp}
 The purpose of this section is to give an additional  perspective on the class of network satisfaction problems of
 finite symmetric integral representable relation algebras  with a flexible atom. 
{Even more, we define these computational problems in this section completely without the use of the relation algebra framework.} 
Our classification  result for these problems does not depend on the content of this section and the reader may skip it.

 We introduce  a class of computational decision problems which we call \emph{network consistency problems (NCPs)}. It is easy to see that NCPs are in a 1-to-1 correspondence with NSPs  of finite, symmetric, integral  $\mathbf{A}\in \rra$ with a flexible atom.

%
%
%

\begin{definition}\label{def:stenc} Let $A$ be a finite set and $R\subseteq A^3$. Then $R$ is called \emph{totally symmetric} if for all bijections $\pi \colon \{1,2,3\}\rightarrow \{1,2,3\}$  we have
	$$(a_1,a_2,a_3)\in R ~\Rightarrow~ (a_{\pi(1)},a_{\pi(2)},a_{\pi(3)})\in R.$$
We call an element  $p\in A$ \emph{identity element} if for all $x,y\in A$ the following holds:$$(p,x,y)\in {R} \Leftrightarrow  x=y.$$
A structure $(A;R)$ is called a \emph{stencil} if $R$  is totally symmetric and it contains an identity element.

	\end{definition}

\begin{definition}Let $(G;F)$ be an undirected graph and let $Q$ be a set. We call a map $c\colon F\rightarrow Q$ an \emph{edge $Q$-coloring of $(G;F)$} if for all $x,y\in G$ with $(x,y)\in F$ it holds that $c((x,y))=c((y,x))$.
	\end{definition}
For each fixed stencil, we define  an NCP as follows.

\begin{definition}\label{def:NCP} Let  $(A,R)$ be a stencil. The \emph{network completion problem of  $(A,R)$},  denoted by $\ncp(A,R)$, is the following problem. Given a finite undirected graph $(G;F)$ with an edge $\mathcal{P}(A)$-coloring $f$ the task is to decide whether there exists an  edge $A$-coloring $f'$ of $(G;F)$  such that
	\begin{enumerate}
		\item for all $x,y\in G$ with $(x,y)\in F$ it holds that $f'((x,y))\in f((x,y))$.
		\item for all $x,y,z\in G$ with $(x,y), (y,z), (x,z)\in F$ we have $$(f'((x,y)), f'((y,z)), f'((x,z)) )\in R.$$
	\end{enumerate}
	\end{definition}

The following proposition illustrates how NCPs correspond to a certain class of NSPs.

\begin{proposition}\label{prop: ncp nsp}  The class of NCPs and 
	the class of NSPs for finite symmetric integral representable relation algebras with a flexible atom are in a natural 1-to-1 correspondence such that corresponding problems are polynomial-time equivalent.
	\end{proposition}

\begin{proof}Let $(A', R')$ be a stencil with $p\in A'$ according to (2) in Definition \ref{def:stenc} and let $s$ be new element with $s\not\in A'$. We define a relational structure $\mathfrak{D}$ as follows. The domain of $\mathfrak{D}$ is the set $A'\cup \{s\}$. We assume that $\mathfrak{D}$ has every subset of its domain as a unary relation. Furthermore $\mathfrak{D}$ contains the binary relation $E^\mathfrak{D}:=\{(x,x)\mid x\in A'\cup\{ s \}\}$  and the ternary relation $H^\mathfrak{D}$ which is defined as follows: 
	\begin{align*}
		(x,y,z)\in H^\mathfrak{D}:\Leftrightarrow & ~\bigl((x,y,z)\in R' \\ & ~\vee (s\in \{x,y,z\}\wedge p\not \in\{x,y,z\}) \\ &~\vee ((x,y,z)\in \{(p,s,s), (s,p,s), (s,s,p)\} ) \bigr).
	\end{align*}
One can find a finite symmetric integral algebra $\mathbf{A}\in \rra$ with domain $\mathcal{P}(A'\cup\{s\})$ and a flexible atom $s$ such that $\mathfrak{D}$ is the atom structure of $\mathbf{A}$ (see, e.g., Theorem~2.2 and Theorem~2.6 in the article by \cite{Maddux1982}).

 Furthermore, given a finite symmetric integral algebra $\mathbf{A}\in \rra$ with a flexible atom $s$ let $\mathfrak{D}$ be the atom structure of $\mathbf{A}$. Let $(A'; R')$ be the substructure of the $\{R\}$-reduct of $\mathfrak{D}$ induced by $D\setminus\{s\}$. By the properties of $\mathfrak{D}$ we get that $(A'; R')$  is a stencil.

We show that the instances of $\ncp(A', R')$ and $\nsp(\mathbf{A})$ are in a natural 1-to-1 correspondence that preserves the acceptance condition of the computational problems.
Let  $(G;F)$ be a finite undirected graph with an edge $\mathcal{P}(A')$-coloring $f$. We define an $\mathbf{A}$-network $(G;g)$ by defining
$$g(x,y) = \left\{
\begin{array}{ll}
	f(x,y) & ~~(x,y)\in F , \\
	\{p\}&~~(x,y)\not\in F \text{~and~} x=y\\
	\{s\} & ~~\, \textrm{else.} \\
\end{array}
\right. $$
It is easy to see that  $(G;F)$ is an accepted instance of  $\ncp(A', R')$  if and only if  $(G;g)$ is an accepted instance of $\nsp(\mathbf{A})$.  Since we can reverse this and find for every  $\mathbf{A}$-network $(G;g)$  a finite undirected graph  $(G;F)$ with an edge $\mathcal{P}(A')$-coloring $f$ such that each of them is an accepted instance  if and only if the other one is. 
These to reductions show that the computational decision problems 
$\ncp(A', R')$ and $\nsp(\mathbf{A})$  are polynomial-time equivalent.
 
%
%
	\end{proof}



	We end this section by providing  a rich source of examples for $\ncp$s.
	
	\begin{example}[``Distance problems'']
		Let $A\subset \mathbb{Q}$ be an arbitrary finite set that contains $0$. We define the relation $R\subseteq A^3$ of all tuples which satisfy all instantiations of the triangle inequality, i.e.
		$$(a_1,a_2,a_3)\in R ~:\Leftrightarrow~ (a_1\leq a_2+a_3)  ~\wedge~  (a_2\leq a_1+a_3)~\wedge~ (a_3\leq a_1 +a_2),$$
		where the addition is meant to be the usual one on rational numbers. The relation $R$ is by definition totally symmetric and the element $0$ is an identity element. Therefore, $(A;R)$ is a stencil. 
		
		Now consider a finite undirected graph $(G;F)$  with an edge $\mathcal{P}(A)$-coloring $f$. This can be seen as a labeling of each edge in the graph by a set of  possible (or allowed) distances. 
		The computational task  of $\ncp(A,R)$ is  to decide whether  one can choose for each edge $(x,y)$ one of the possible distance such that in the end this choice satisfies on each triangle of edges the triangle inequalities of metric spaces.

		By Proposition~\ref{prop: ncp nsp} there exists a finite symmetric integral algebra  $\mathbf{A'}\in \rra$ with a flexible atom $s$ such that $\ncp(A,R)$  and $\NSP(\mathbf{A'})$ are polynomial-time equivalent. By the proof of Proposition~\ref{prop: ncp nsp} we have that the domain of $\mathbf{A'}$ is equal to  $\mathcal{P}(A\cup\{s\})$. It is easy to observe that $A\cup \{s\}$ is the set of atoms $A'_0$ from $\mathbf{A'}$. 
		
		We define an operation $f\colon {A'_0}^6\rightarrow A'_0$ as follows:
		
		$$f(x_1,\ldots,x_6) = \left\{
		\begin{array}{ll}
			  s & ~~s \in \{x_1,\ldots, x_6 \}, \\
			\max\{x_1,\ldots,x_6\}&~~\text{otherwise},\\
		\end{array}
		\right. $$where the $\max$ operation is the usual from in $\mathbb{Q}$.
		
		The allowed  triples of $\mathbf{A'}$ are, up to triples that involve the flexible atom $s$, those which arise from valid triangle inequalities. For this reason the operation $f$ preserves the allowed triples of $\mathbf{A'}$. Moreover, one can check that $f$ satisfies the Siggers identity. 
		This implies that all these ``distance problems'' satisfy the first condition in Theorem~\ref{theo:result2}  and are therefore solvable by a polynomial-time algorithm.
		{Furthermore, if the normal representation of  $\bf{A}' $ induces  a proper relation algebra $\bf{A}$ on the set $A$, then $\bf{A}$ is a representable relation algebra with a normal representation.   This follows from 
		 		Proposition 2.7.4 in the thesis of  \cite{Conant2015} (see also~\cite{Del}), since the composition operation of $\bf{A}$ is associative.
		 	We  get that  if  the representable relation algebra $\bf{A}$ exists then the argument from Example~\ref{examp:+fexible atom}  implies that $\nsp(\bf{A})$  is also in P.}
		 
	\end{example}


\section{Binary Injective Polymorphisms}\label{section: can binary injection}
We give in this section a proof of the following proposition.

\begin{proposition}\label{binaryinjectionNPhard}
	Let $\mathfrak{B}$ be a normal representation of a finite, symmetric, integral  $\mathbf{A}\in \rra$ with a flexible atom $s$. If $\HSPf(\{\Pol(\mathfrak{B }) \}  )$ does not contain a 2-element algebra where all operations are projections, then $\mathfrak{B}$ has a binary injective polymorphism. 
\end{proposition}
{This statement is a consequence of well known results that can be found in the book by \cite{Book} and results by \cite{MottetPinskerSmooth} applied to the class of normal representations of finite, symmetric, integral $\mathbf{A}\in \rra$ with a flexible atom $s$.}
An operation $f \in \Pol(\mathfrak{B})$ is called \emph{essentially unary} if it depends on at most one of its variables and  $f$ is called \emph{essential} otherwise.

Following the terminology of \cite{MottetPinskerSmooth}, we now define
\emph{free  2-orbits}. 
{The existence of a free 2-orbit appeared under the name \emph{`orbital extension property'} for example in the habilitation thesis of \cite{Bodirsky-HDR-v8}.}

\begin{definition}
Let $\mathfrak{B}$ be a structure. A 2-orbit $O$  of $\Aut(\mathfrak{B })$ is called \emph{free} if for all elements $x,y\in {B}$ there exists $z\in B$ with $(z,x)\in O $ and $(z,y)\in O$.
\end{definition}

Note that if $\Aut({\mathfrak B})$ has a free 2-orbit then it is transitive. The following theorem generalises a fact that was first proved for first-order reducts of $({\mathbb Q};<)$  by \cite{tcsps-journal}.

\begin{proposition}[Lemma 5.3.10 in~\cite{Bodirsky-HDR-v8}]\label{prop:binary essential}
	Let $\mathfrak{B}$ be a structure such that
	$\mathfrak{B}$ has a free 2-orbit.
	If $\Pol(\mathfrak{B})$ contains an essential operation then it contains a binary essential operation.
	\end{proposition}


The following is essentially taken from the article by \cite{MottetPinskerSmooth}. 

\begin{definition}
Let $\mathfrak{B}$ be a structure.
Then the \emph{canonical binary structure} of
$\mathfrak B$ is the structure with domain $B$ and a binary relation for each 2-orbit $O$ of $\Aut(\mathfrak{B})$ such that $(x,y)\in O$ implies $x\neq y$.
	\end{definition}

\begin{definition}A $\tau$-structure $\mathfrak{B}$ has \emph{finite duality} if there exists a finite set $\mathcal{F}$ of finite $\tau$-structures such that a $\tau$-structure $\mathfrak{I}$ has a homomorphism to $\mathfrak{B}$ if and only if no element of $\mathcal{F}$ has a homomorphism to $\mathfrak{I}$.
	\end{definition}

We establish finite duality for the class of structures which is important for our classification purposes.

\begin{lemma}\label{lem: finite duality}
	Let $\mathfrak{B}$ be a normal representation of a finite, integral  $\mathbf{A}\in \rra$ with a flexible atom $s$.  Then the canonical binary structure of $\mathfrak{B}$ has finite duality.
\end{lemma}
\begin{proof}Let $\tau := A_0\setminus\{\id\}$. 
	Note that since $\mathbf{A}$ is integral and $\mathfrak{B}$ is homogeneous, the canonical binary structure $\mathfrak{C}$ of  $\Aut(\mathfrak{B})$ is precisely the $\tau$-reduct of $\mathfrak{B}$. 
	Let  $\mathcal{F}$ be the set  of all $\tau$-structures with domain $\{1,2,3\}$ that do not have a homomorphism to $\mathfrak{C}$.
We show that the set $\mathcal{F}$ witnesses the  finite duality of the canonical binary structure $\mathfrak{C}$ of  $\mathfrak{B}$.
	Let  $\mathfrak{I}$ be a $\tau$-structure with a homomorphism to $\mathfrak{C}$.
 If there exists $\mathfrak{F}\in \mathcal{F}$ with a homomorphism to $\mathfrak{I}$, then $\mathfrak{F}$ also has  a homomorphism to $\mathfrak{C}$, contradicting the choice of $\mathcal{F}$.
	For the other direction assume that no element from $\mathcal{F}$ has a homomorphism to $\mathfrak{I}$. Let $\mathfrak{I}'$ be the $\tau$-expansion of the $(\tau \setminus\{s\})$-reduct of $\mathfrak{I}$ where the relation $s^{\mathfrak{I}'}$ is defined by $$s^{\mathfrak{I}'}:=\{(x,y)\in I^2 \mid  (x,y)\in s^{\mathfrak{I}} \vee (x\neq y \wedge \forall a\in \tau \setminus\{s\}. \, (x,y) \notin a^\mathfrak{I} )\}.$$
	By the definition of the flexible atom $s$ it follows that no element from $\mathcal{F}$ has a homomorphism to $\mathfrak{I}'$. This implies  that for all distinct elements $x,y$ from  $\mathfrak{I}'$ the tuple $(x,y)$ is in at most one relation from $\tau$.
	The definition of $\mathfrak{I}'$ ensures that $(x,y)$ is in at least one relation from $\tau$.
	 Recall that  Proposition \ref{prop: finite boundedness} states that $\mathfrak{C}$ is finitely bounded by $\tau$-structures of size at most three. Assume that one of those bounds $\mathfrak{N}$ embeds into $\mathfrak{I}'$. This implies by what we noted before that all elements of  $\mathfrak{N}$ are in precisely one relation from $\tau$. On the other hand $\mathfrak{N}$ is not in $\mathcal{F}$ and therefore has  a homomorphism to $\mathfrak{C}$. Since all elements of $\mathfrak{N}$ are related by precisely one relation from $\tau$, this homomorphism needs to be a embedding,  contradicting our  assumption on $\mathfrak{N}$ to be a bound. Therefore, none of the bounds of $\mathfrak{C}$ embeds into $\mathfrak{I}'$, which means that $\mathfrak{I}'$ is a substructure of $\mathfrak{C}$. Clearly, there exists a homomorphism from $\mathfrak{I}$ to $\mathfrak{I}'$ which proves the lemma. \end{proof}

The following proposition about the existence of injective operations is from \cite{MottetPinskerSmooth}, building on ideas of  \cite{BodPin-Schaefer-both} and \cite{BMPP16}.
\begin{proposition}[\hspace{1sp}\cite{MottetPinskerSmooth}]\label{prop:approxinjective}
	Let $\mathfrak{B}$ be a homogeneous structure such that $\Aut(\mathfrak{B})$ is transitive and such that the canonical binary structure of  $\Aut(\mathfrak{B})$ has finite duality.
	If $\Pol(\mathfrak{B})$ contains a binary essential operation that preserves $\neq$ then it contains a binary injective operation.
	\end{proposition}

We are now able to prove the main result of this section.
\begin{proof}[Proof of Proposition \ref{binaryinjectionNPhard}]
Note that since $\mathbf{A}$ is integral and $\mathfrak{B}$ is homogeneous, the flexible atom $s$ is a free 2-orbit of $\Aut(\mathfrak{B})$. Furthermore, $\Aut(\mathfrak{B})$ is transitive. 
	Suppose that $\HSPf(\{\Pol(\mathfrak{B }) \}  )$ does not contain a 2-element algebra where all operations are projections. 
	Since all operations of $\Pol(\mathfrak{B})$ are edge conservative, it follows that 
	$\Pol(\mathfrak B)$ contains an operation that does not behave as a projection on $\{s,\id\}$. 
	This implies that 	$\Pol(\mathfrak B)$ contains an essential operation. 
	By Proposition~\ref{prop:binary essential},
	$\Pol(\mathfrak B)$ must also contain a binary essential operation. Since the canonical binary structure of $\mathfrak{B}$ has  finite duality by Lemma~\ref{lem: finite duality} we can apply Proposition~\ref{prop:approxinjective} and get that $\Pol(\mathfrak{B})$  contains a binary injective operation. 
 \end{proof}

The following shows how  to use Proposition~\ref{binaryinjectionNPhard} to obtain a hardness result  for the concrete  $\mathbf{A}\in \rra$ from Example \ref{exam:henson}.

\begin{example}[Hardness of representable relation algebra  $\#17$, see~\cite{BMPP16,BodirskyKnaeRamics}] 
	Let $\mathfrak{N}'$ be the normal representation of the algebra $\#17$ mentioned in Example \ref{exam:henson}. We claim that the structure $\mathfrak{N}'$ does not have a binary injective polymorphism. To see this, consider a substructure of $\mathfrak{N}'^2$ on elements $x,y,z\in V^2$ such that $ (E,=)(x,y)$, $(= ,E)(y,x)$, and $(E,E)(x,z)  $. Assume  $\mathfrak{N}'$ has a binary injective  polymorphism $f$.  This means that $\overline{f}(E,\id)=E=\overline{f}(\id,E)$ holds. Then we get that $E(f(x), f(y))$, $E(f(y),f(z))$, and $E(f(x),f(z)$ hold in $\mathfrak{N}'$, which is a contradiction, since in $\mathfrak{N}'$ triangles of this form are forbidden. 
	By the contraposition of  Proposition~\ref{binaryinjectionNPhard}	it follows that
	 $\HSPf(\{\Pol(\mathfrak{B }) \}  )$  contains a 2-element algebra where all operations are projections.
We conclude with Theorem~\ref{theo:hardness} that  $\nsp(\#17)$ is an NP-hard problem.	
	
\end{example}

\section{From Partial to Total Canonical Behaviour}\label{sec:canon}
 In this section we prove that in many cases the existence of a polymorphism with a certain partial behaviour implies the existence of a canonical 
	 polymorphism with the same partial behaviour. 
	Following this idea we start in Section~\ref{subsec:1} with the proof that the existence of  an injective polymorphism implies the existence of a canonical injective polymorphism.
In some cases the existence of an $\{a,b\}$-canonical polymorphism implies the existence of a canonical polymorphism with the same behaviour on $\{a,b\}$.
We prove this separately for binary (Section~\ref{subsec:2}) and ternary (Section~\ref{subsec: 3})  operations, making use of the binary injective polymorphism that exists by the results from Section~\ref{section: can binary injection} and Section~\ref{subsec:1}.

Let us remark that most proofs of this section would fail if the representable relation algebra $\mathbf{A}$ was not symmetric.  Indeed, every representable relation algebra that contains a non-symmetric atom a normal representation  would not satisfy Proposition \ref{prop:binary poly from order to B} below, since the stated behaviour  is not well-defined.

We assume for this section that $\mathfrak{B}$ is a normal representation of a finite, symmetric, integral $\mathbf{A}\in \rra$ with a flexible atom $s$. Let furthermore $\mathfrak{B}_<$ be the expansion of $\mathfrak{B }$ by the generic linear order. The structure $\mathfrak{B}_<$ exists by the observations in Section \ref{subsec: normalre}. 




\begin{proposition}\label{restriction} 
		Let $f \in \Pol(\mathfrak{ B})^{(n)}$ be injective.  Then there exists a polymorphism $f_<$ of $\mathfrak{B}_<$ and an injective endomorphism $e$ of $\mathfrak{B}$ such that $$ f=e\circ f_<$$ 	as mappings from $B^n$ to $B$.
\end{proposition}

\begin{proof}

Let $U:=f(B^n)$ and consider  the substructure $\mathfrak{U}$ induced by  $\mathfrak{B}$ on $U$.
There exists a linear ordering on $B^n$, namely the lexicographic order given by the linear order of $\mathfrak{B}_<$ on each coordinate.

 Let $\mathfrak{U}_<$ be the expansion  of $\mathfrak{U}$ by the linear order that is induced by the lexicographic linear order  of $\mathfrak{B}_<$ on the preimage. This is well defined since $f$ is injective. By the definition of $\mathfrak{B}_<$ and a compactness argument the structure $\mathfrak{U}_<$ embeds into $\mathfrak{B}_<$. In this way we obtain a homomorphism $f_<$  from $\mathfrak{B}_<^n$  to $\mathfrak{B}_<$. Again by a compactness argument also an endomorphism $e$ with the desired properties exists. \end{proof}

 \subsection{Canonical Binary Injective Polymorphisms}\label{subsec:1}

%


We prove in this section that the existence of an injective polymorphism implies the existence of a canonical injective polymorphism.
We say that a polymorphism $f$ of $\mathfrak{B}_<$ is \emph{canonical with respect to $\mathfrak{B}_<$} if $f$ satisfies Definition \ref{defi:canonicalfunctions}, where the underlying representable relation algebra is the proper relation algebra induced by the binary first-order definable relations (i.e., unions of  2-orbits) of $\mathfrak{B }_<$.
Note that in a normal representation  the set of $2$-orbits equals the set of the interpretations of the atoms of the representable relation algebra.

\begin{proposition}\label{prop:binary poly from order to B}
Let $f$ be a binary polymorphism of $\mathfrak{B}_<$ that is canonical with respect to $\mathfrak{B}_<$. 
{Let  $h \colon A_0^2 \rightarrow A_0$ be the map such that for all $x,y,z\in A_0$
$$ h(x,y)=z  \Leftrightarrow \bar{f}(x^{\mathfrak{B}_<} \cap \leq^{\mathfrak{B}_<}, y^{\mathfrak{B}_<} \cap \leq^{\mathfrak{B}_<}) =z^{\mathfrak{B}_<} \cap \leq^{\mathfrak{B}_<}, $$ 
(cf. Remark \ref{rem:proper RA Bor}).    
Then $h$ is well defined and there exists a canonical binary polymorphism of $\mathfrak{B}$ with behaviour $h$.}

\end{proposition}
\begin{proof}
	
	The function $h$ is well defined since all atoms are symmetric. We show that there exists a canonical polymorphism of $\mathfrak{B}$ that has $h$ as a behaviour. 
	Consider the following structure $\mathfrak{A}$ on the domain $B^2$. Let $x,y\in B^2$ and let $a,a_1,a_2\in A_0$ be atoms of $\mathbf{A}$ with $a_1^\mathfrak{B}(x_1,y_1)$ and $a_2^\mathfrak{B}(x_2,y_2)$. Then we define that $a^\mathfrak{A}(x,y)$ holds  if and only if $h(a_1,a_2)=a$.

	We show in the following that $\mathfrak{A}$ has a homomorphism to $\mathfrak{B}$. This is enough to prove the statement, because a homomorphism from $\mathfrak{A}$ to $\mathfrak{B}$ is a canonical polymorphism of $\mathfrak{B}$. 
	Since $\mathfrak{B}$ is homogeneous  is suffices to show that every finite substructure of $\mathfrak{A}$ homomorphically maps to $\mathfrak{B}$.
	Let $\mathfrak{F}$ be a finite substructure of $\mathfrak{A}$ and assume for contradiction that $\mathfrak{F}$ does not homomorphically map to $\mathfrak{B}$. We can view $\mathfrak{F}$ as an atomic $\mathbf{A}$-network. Since $\mathfrak{B}$ is fully universal $\mathfrak{F}$ is not  closed. There must exist elements $ b^1,b^2, b^3 \in B^2 $ of $\mathfrak{F}$ and atoms $a_1,a_2,a_3 \in A_0$ such that $a_1  \not\leq a_2\circ a_3$ holds in $\mathbf{A}$ and
	$$a_1^\mathfrak{F}(b^1,b^3), ~ a_2^\mathfrak{F}(b^1,b^2), \text{~and~} a_3^\mathfrak{F}(b^2,b^3). $$
	This means that the substructure induced on the elements  $ b^1,b^2, b^3 $ by $\mathfrak{F}$ contains a forbidden triple. 
	
	Now we consider the substructures that are induced on $ b_1^1,b_1^2, b_1^3 $ and  $ b_2^1,b_2^2, b_2^3 $  by $\mathfrak{B}$. Our goal is to order these elements such that for all $i,j \in \{1,2,3\}$ 
	\begin{equation}\label{equation}
\neg (b_1^i< b_1^j \wedge  b_2^i > b_2^j).
	\end{equation}
	If we achieve this we know that there exist elements in $\mathfrak{B}_<$ that induce isomorphic copies of the induced structures of the elements  $ b_1^1,b_1^2, b_1^3 $ and  $ b_2^1,b_2^2, b_2^3 $ with the additional ordering.
	Now the application of the polymorphism $f$ on these elements results in a structure whose $A_0$-reduct  is isomorphic to the substructure induced by $b^1,b^2$ and $b^3$ on $\mathfrak{F}$ by the definition of the canonical behaviour $h$. This contradicts our assumption because a polymorphism can not have a forbidden substructure in its image.
	
	It remains to show that we can choose orderings on the elements  $ b_1^1,b_1^2, b_1^3 $ and  $ b_2^1,b_2^2, b_2^3$ such that (\ref{equation}) holds.	Without loss of generality we can assume that $\{b_1^1,b_1^2, b_1^3 \}\cap \{b_2^1,b_2^2, b_2^3\} = \emptyset$ holds. Now consider the following cases:
	
	\begin{enumerate}
		\item $|\{b_1^1,b_1^2, b_1^3 \}|=3$ and $|\{b_2^1,b_2^2, b_2^3\} |= 3$.
		
		We can obviously choose  linear orders on both sets such that (\ref{equation}) holds.
		\item $|\{b_1^1,b_1^2, b_1^3 \}|=2$ and $|\{b_2^1,b_2^2, b_2^3\} |= 3$.
		
		Assume that $\id^\mathfrak{B}(b_1^1,b_1^2)$ holds then the possible orders are
		$$ b_1^1=b_1^2 < b_1^3 \textup{~and~} b_2^1<b_2^2 <b_2^3.$$
		\item $|\{b_1^1,b_1^2, b_1^3 \}|=2$ and $|\{b_2^1,b_2^2, b_2^3\} |= 2$.
		
		
	First consider the case that $\id^\mathfrak{B}(b_1^1,b_1^2)$ and $\id^\mathfrak{B}(b_2^1,b_2^2)$ hold. Then we choose as orders 
			$$ b_1^1=b_1^2 < b_1^3 \textup{~and~} b_2^1=b_2^2 <b_2^3.$$
			
			In the second possible case we can assume without loss of generality that
			$\id^\mathfrak{B}(b_1^1,b_1^2)$ and $\id^\mathfrak{B}(b_2^2,b_2^3)$ hold. Note that otherwise we could change the role of  two of the tuples $b^1,b^2$ and $b^3$ and get this case. The compatible order is then
				$$ b_1^1=b_1^2 < b_1^3 \textup{~and~} b_2^1<b_2^2 =b_2^3.$$
				\item $|\{b_1^1,b_1^2, b_1^3 \}|=1$ and $|\{b_2^1,b_2^2, b_2^3\} |= 3$.\\
		In this case we choose the order
			$$ b_1^1=b_1^2 = b_1^3 \textup{~and~} b_2^1<b_2^2 <b_2^3.$$
			\item 	 $|\{b_1^1,b_1^2, b_1^3 \}|=1$ and $|\{b_2^1,b_2^2, b_2^3\} |= 2$.\\
			Assume that $\id^\mathfrak{B}(b_2^1,b_2^2)$ holds and that we have
				$$ b_1^1=b_1^2 = b_1^3 \textup{~and~} b_2^1=b_2^2 <b_2^3.$$
				\item $|\{b_1^1,b_1^2, b_1^3 \}|=1$ and $|\{b_2^1,b_2^2, b_2^3\} |= 1$.\\
				For this case we trivially get 
					$$ b_1^1=b_1^2 = b_1^3 \textup{~and~} b_2^1=b_2^2 =b_2^3.$$
	\end{enumerate}
	Note that up to the symmetry of the arguments for both coordinates these are all the possible cases. This completes the proof of the proposition.
 \end{proof}


\begin{corollary}\label{injectivecanonical poly}
	Suppose that $\mathfrak{B}$ has a binary injective polymorphism. Then $\mathfrak{B}$ also has a canonical binary injective polymorphism.
\end{corollary}
\begin{proof}
By Proposition \ref{restriction} we may assume that there exists also an injective polymorphism of $\mathfrak{B }_<$. 
	The structure $\mathfrak{B }_<$ has the Ramsey property by Theorem~\ref{theorem:ramsey order}.
	Therefore,  Theorem \ref{canonisation} implies that there also exists an injective canonical polymorphism $g$ of $\mathfrak{B }_<$. According to Proposition \ref{prop:binary poly from order to B} the restriction of the behaviour $\bar{g}$ 
	to the 2-orbits that satisfy $x\leq y$ induces the behaviour of a canonical polymorphism of $\mathfrak{B }$ which is also injective. \end{proof}

\subsection{Canonical $\{a,b\}$-symmetric Polymorphisms}\label{subsec:2}
 We will now use the results about binary injective polymorphism from Section~\ref{subsec:1} to show the existence of a canonical  $\{a,b\}$-symmetric polymorphism in case there exists an $\{a,b\}$-symmetric polymorphism. 

\begin{lemma}\label{lem: ab symm implies injective}	Let $a,b\in A_0\setminus \{\id\}$ be atoms. Then every binary $\{a,b\}$-symmetric polymorphism  of $\mathfrak{B}$ is injective.
\end{lemma}

\begin{proof}{ 
	Let $f$ be an $\{a,b\}$-symmetric polymorphism. Without loss of generality  $\bar{f}(a,b)=a=\bar{f}(b,a)$. 
		Assume for contradiction that $f$ is not injective. This means that there exist $c\in A_0$ and $x,y\in B^2$  with $(c,\id)(x,y)$ {(for the notation see Definition~\ref{defi: 2n relation})} such that $\id(f(x),f(y))$ holds. 
		
	\underline{Case 1:} $s \not \in \{a,b\}$. 
		Since $s$ is a flexible atom we may choose $z\in B^2$ such that $(a,b)(z,x)$ and $(s,b)(z,y)$ hold. 
		By the choice of the polymorphism $f$ we get $a(f(z),f(x))$ and $(s\cup b)(f(z),f(y))$ which induces either the forbidden triple $(\id,s,a)$ or the forbidden triple $(\id,b,a)$ on $f(x), f(y) $, and $f(z)$. 
		
\underline{Case 2:} $s=a$. 
 We choose  $z\in B^2$ such that $(a,b)(z,x)$ and $(b,b)(z,y)$. 
This is possible since $a$ is the flexible atom. We obtain  $a(f(z),f(x))$ and 
$ b(f(z),f(y))$ which  again induces a forbidden triple on $f(x), f(y) $, and $f(z)$.

\underline{Case 3:} $s=b$. 
 We choose  $z\in B^2$ such that $(a,b)(z,x)$ and $(b,b)(z,y)$. 
This is possible since $a$ is the flexible atom. We obtain  $a(f(z),f(x))$ and 
$ b(f(z),f(y))$ which again  induces a forbidden triple on $f(x), f(y) $, and $f(z)$.

Since we obtained in all cases a contradiction we conclude that $f$ is injective.}
 \end{proof}

\begin{proposition}\label{partiallysymm}
	Let $a,b\in A_0\setminus \{\id\}$ be atoms.
If $\mathfrak{B}$ has  a binary $\{a,b\}$-symmetric polymorphism, then $\mathfrak{B }$ has  also a binary canonical $\{a,b\}$-symmetric polymorphism.
\end{proposition}
\begin{proof}{
Let $f$ be the binary $\{a,b\}$-symmetric polymorphism. By Lemma \ref{lem: ab symm implies injective} we know that $f$ is injective. By Proposition~\ref{restriction} it induces  a polymorphism $f_<$  on $\mathfrak{B }_<$.
The structure $\mathfrak{B }_<$ has the Ramsey property by Theorem~\ref{theorem:ramsey order}.
  Let $g$ be the canonization of $f_<$  that exists by Theorem \ref{canonisation}. 
The restriction of the behaviour $\bar{g}$ 
to the 2-orbits that satisfy $x\leq y$ induces  by Proposition \ref{prop:binary poly from order to B} the behaviour of a canonical polymorphism $h$ of $\mathfrak{B }$. The way we obtained $h$ ensures that $h$ is $\{a,b\}$-symmetric with the same behaviour on $\{a,b\}$ as $f$. } \end{proof}

{The following is an easy observation about $\{a,\id \}$-symmetric polymorphisms that we will use several times.}

\begin{observation}\label{lemma:a id symm means a}
	Let $a \not\leq \id$ be an atom and $f$ an $\{a,\id \}$-symmetric polymorphism of $\mathfrak{B}$.
	 Then $\bar{f}(a,\id)=a=\bar{f}(\id,a)$.
\end{observation}
\begin{proof}Suppose for contradiction that $\bar{f}(a,\id)=\id=\bar{f}(\id,a)$.
Let  $x,y,z \in B^2$  be such that $$(a,\id)(x,y), ~ (\id,a)(y,z), \text{~ and ~} (a,a)(x,z)$$ and consider the substructure of $\mathfrak{B}$  that is induced by  $f(x), f(y)$ and $f(z)$. This structure induces a forbidden triple $(\id,\id,a)$ which contradicts the assumption $\mathbf{A}\in \rra.$ \end{proof}
	In the remainder of the section we combine canonical $\{a,b\}$-symmetric polymorphisms to obtain a single ``maximal-symmetric'' polymorphism.

\begin{definition}
	We call a subset $\{a,b\} \subseteq A_0$ an \emph{edge of $\Polc(\mathfrak{B })$}
	and we call the elements in  $$Q:= \big\{ \{a,b \} \subseteq A_0  \mid\exists g\in \Polc(\mathfrak{B }) \text{~such that~} \bar{g} \text{~is symmetric on~} \{a,b\}   \big\}$$
	the \emph{red edges} of $\Polc(\mathfrak{B })$. 
	
\end{definition}
The terminology of the colored edges as well as the following lemma are  from \cite{Bulatov-Conservative-Revisited}.

\begin{lemma}\label{maximal-symmetric defi}
	There exists a binary canonical polymorphism 
	 that is symmetric on all red edges and behaves on each non-red edge like a projection. We call this function \emph{maximal-symmetric}.
\end{lemma}
\begin{proof}
	For each $\{a,b\} \in Q$ let $f_{a,b}$ be a canonical polymorphism such that its behaviour  is symmetric on $\{a,b\}$.
	We prove the lemma by an induction  on the size of subsets of $Q$, i.e., we show that for every subset $F\subseteq Q$ of size $n$ there exists a polymorphism $f_F\in \Polc(\mathfrak{B })$ that is symmetric on all edges from $F$.
	For each subset $\{a\}$ of $Q$  of size one, there exists by the definition of red edges a canonical  polymorphism $f_{a,a}$ with a behaviour that is symmetric on $\{a\}$. Let $F\subseteq Q$ and suppose there exists a canonical polymorphism $g$ with symmetric behaviour on elements from $F$. Let $\{a_1,a_2\} \in Q\setminus F$.  We want to show that there exists a canonical polymorphism with a behaviour that is symmetric on all elements from $F\cup \{a_1,a_2\}$.
	We may assume that this does not hold for $g$, otherwise we are done. Therefore, and since $g$ is  edge-conservative, we have$$\bar{g}(a_1,a_2) \not = \bar{g}(a_2,a_1) \text{~and ~} \bar{g}(a_1,a_2),~ \bar{g}(a_2,a_1) \in \{a_1,a_2\}.$$
	With this it is easy to see that  $f_{a_1,a_2}(g(x,y),g(y,x ))$ is a polymorphism with a behaviour that is symmetric on all elements from $F\cup \{a_1,a_2\}$.

	This proves the first part of the statement. For the second part note that for a binary canonical edge-conservative polymorphism there are only 4 possibilities for the behaviour on a set $\{a_1,a_2\}$. If $\{a_1,a_2\}$ is not a red edge then every binary canonical  edge-conservative polymorphism behaves like a projection on $\{a_1,a_2\}$ .
	 \end{proof}

\subsection{Canonical Ternary Polymorphisms}\label{subsec: 3}
We obtain in this section a result  that states that the existence of a ternary $\{a,b\}$-canonical polymorphism $f$ implies the existence of a canonical polymorphism with the same behaviour on $\{a,b\}$ as $f$ (Corollary~\ref{partialto global:majomino}). This is done similarly as in   Section~\ref{subsec:2}.

\begin{lemma}\label{s lemma}
	Let $s'\in \Polc(\mathfrak{B })$ be a binary, injective, maximal-symmetric polymorphism. Then the function $s^*\colon  \mathfrak{B}^3 \rightarrow\mathfrak{B}^3$ where $s^*(x_1,x_2,x_3) $ is defined by
	$$ ( s'(s'(x_1,x_2),s'(x_2,x_3) ), s'(s'(x_2,x_3),s'(x_3,x_1)),s'(s'(x_3,x_1),s'(x_1,x_2))     )$$
	is a homomorphism. 
	Moreover, for all  $x,y\in B^3$ with $x\not =y$ it holds that $$\overline{\id}^{\mathfrak{B}^3}(s^*(x),s^*(y)).$$
\end{lemma}
Note that this means that two distinct tuples in the image of $s^*$ have distinct entries in each coordinate.

\begin{proof}
	Let $x,y \in B^3$ and 
	suppose  that $a^{\mathfrak{B}^3}(x,y)$ holds for $a\in A$. By the definition of the product structure $a^\mathfrak{B}(x_i,y_i)$ holds for all $i\in\{1,2,3\}$.
	Since $s'$ is a polymorphism clearly $a^\mathfrak{B}(s^*(x)_i, s^*(y)_i     )$ holds by the definition of $s^*$. Now we use again the definition of a product structure and get $a^{\mathfrak{B}^3}(s^*(x), s^*(y)    )$ which shows that $s^*$ is a homomorphism.

	For the second part of the statement let  $x,y \in B^3$ distinct. Suppose that $a_1^\mathfrak{B}(x_1,y_1) $, $a_2^\mathfrak{B}(x_2,y_2)$ and $a_3^\mathfrak{B}(x_3,y_3)$ hold for some $a_1,a_2,a_3\in A_0$, where at least one atom is different from $\id$.
	Since $s'$ is injective we have 
	$$\overline{\id}^\mathfrak{B}({s'}(x_1,x_2) ,{s'}(y_1,y_2)) \text{~or~ }{\overline{\id}}^\mathfrak{B}({s'}(x_2,x_3) ,{s'}(y_2,y_3)).
	$$
Again by the injectivity of $s'$ we get $$\overline{\id}^\mathfrak{B}(   s'(  {s'}(x_1,x_2) , {s'}(x_2,x_3)    )       ,    s'( {s'}(y_1,y_2), s'(y_2,y_3))).$$
	By the definition of $s^*$ this shows that $\overline{\id}^\mathfrak{B}(s^*(x)_1, s^*(y)_1 )$ holds. It is easy	to see by analogous arguments that the same is true for the other coordinates. Therefore, the statement follows. \end{proof}
	
	

\begin{proposition}\label{prop:partical to total ternery 1.}
	Let $a\not \leq \id$ and $b\not \leq \id$ be atoms of $\mathbf{A}$ such that $\{a,b\} \not \in Q$. Let $m \in \Pol(\frak B)$ be ternary $\{a,b\} $-canonical and $s'\in \Pol(\frak B)$ be injective, maximal-symmetric. Then there exists a canonical $m'\in \Pol(\frak B)$ with the same behaviour as $m$ on $\{a,b\} $.
\end{proposition}

\begin{proof}
By Lemma~\ref{maximal-symmetric defi} we may assume that $s'$ behaves on $\{a,b\}$ like the projection to the first coordinate since $\{a,b\} \not \in Q$. Let $s^*$ be the function defined in Lemma \ref{s lemma} and consider the function $m'\colon B^3\rightarrow B$ which is defined by $m^*(x):=m(s^*(x)) $.

		\medskip
	
	\underline{Claim 1:} $m^*$ is injective.
Let $x,y \in B^3$ be two distinct elements. By Lemma~\ref{s lemma} we know that 
$\overline{\id}^{\mathfrak{B}^3}(s^*(x),s^*(y))$ holds.
Since $m$ is a polymorphism of $\mathfrak{B}$ we directly get that $\overline{\id}^{\mathfrak{B}}(m^*(x),m^*(y))$ holds, which proves the injectivity of $m^*$.
	\medskip
	
	\underline{Claim 2:} $m^*$ is $\{a,b\}$-canonical and behaves on $\{a,b\}$ like $m$.
 
	Let $x,y \in B^3$  with $(q_1,q_2,q_3)(x,y)$ such $q_1,q_2,q_3 \in \{a,b\}$.
	Since $s$ behaves like the first projection on $\{a,b\}$  it follows  that $(q_1,q_2,q_3)(s^*(x),s^*(y))$.
Together with the $\{a,b\}$-canonicity of $m$ this proves Claim 2.

	\medskip
 Since  $m^*$ is injective there exists by Proposition \ref{restriction} a  polymorphism $m^*_<$ of $\mathfrak{B}_<$. 
 Since $\mathfrak{B}_<$ is a Ramsey structure we can apply Theorem~\ref{canonisation} to $m^*_<$. Let $g$ be the resulting polymorphism that is canonical with respect to $\mathfrak{B}_<$. Note that if we consider  $g$ as a polymorphism of $\mathfrak{B}$ it behaves on $\{a,b\} $ like  $m^*$ and therefore like $m$.
 Now  we consider the induced behaviour of $g$ on all 2-orbits that satisfy $x <y$.
Since all atoms of $\mathbf{A}$ are symmetric and $\bar{g}$ is conservative this induces a function $h\colon (A_0\setminus\{\id\}) ^3 \rightarrow A_0\setminus\{\id\}$.
		\medskip
		
	\underline{Claim 3:} The partial behaviour $h$ does not induce a forbidden triple.

We would like to note that we prove this claim with similar arguments as in the proof of Proposition~\ref{prop:binary poly from order to B}.	Assume for contradiction that there exist $x,y,z \in B^3$ such that the application of an operation with behaviour $h$ would induce a forbidden triple.
	Without loss of generality we can order the elements of each coordinate of $x,y,z $ strictly with $x_i<y_i<z_i$ for $i\in \{1,2,3\}$. Note that if on some coordinate there would be the relation $\id$ then we are out of the domain of the behaviour $h$.
	
	If we choose such an order we can find isomorphic copies $\mathfrak{A}$ of this structure (with the order) in $\mathfrak{B}_<$. If we apply the polymorphism $g$ to this copy and forget the order of the structure $g(\mathfrak{A})$ we get a structure that is by definition isomorphic to the forbidden triple.
	This proves Claim 3.
	
		\medskip
	To finish the proof of the lemma note that the composition of $s^*$ with the projection to the $i$-th coordinate for $i\in \{1,2,3\}$ is a canonical, injective polymorphism (for injectivity see Lemma~\ref{s lemma}) and therefore induces a behaviour $f_i\colon A_0^3\rightarrow  A_0\setminus\{\id\}$.
	We define   $f\colon A_0^3\rightarrow  (A_0\setminus\{\id\}) ^3$  by $f(a_1,a_2,a_3):=(f_1(a_1),f_2(a_2),f_3(a_3))$. 
	The composition $h\circ f \colon A_0^3 \rightarrow A_0$ is the behaviour of a canonical function of $\mathfrak{B}$. If $h\circ f $ would induce a forbidden triple 
	then also $h$ would induce a forbidden triple, 
	which contradicts Claim 3. \end{proof}

\begin{corollary}\label{partialto global:majomino}
{Let $\frak B$ have  a binary  injective  polymorphism.}
	Let $a,b\in A_0$ be such that no $\{a,b\} $-symmetric polymorphism exists.  Let $m$ be a ternary $\{a,b\} $-canonical polymorphism. Then there exists a canonical polymorphism $m'$ with the same behaviour on $\{a,b\} $ as $m$.
\end{corollary}
{
\begin{proof}
By assumption there is no canonical $\{a,b\} $-symmetric polymorphism and therefore  $\{a,b\} \not \in Q$. By Corollary \ref{injectivecanonical poly} $\mathfrak{B }$ has a canonical binary injective  polymorphism. This polymorphism is a witness that $\{c,\id\}\in Q$  for all $c\in A_0\setminus \{\id\}$.
	With Lemma \ref{maximal-symmetric defi} we get a maximal symmetric polymorphism $h$. Since $\{c,\id\}\in Q$ we get that $h$ is $\{c,\id\}$-symmetric for all $c\in A_0$. By Observation~\ref{lemma:a id symm means a} it follows $\bar{h}(c,\id)=c=\bar{h}(\id,c)$, which implies that $h$ is injective.
Now  Proposition \ref{prop:partical to total ternery 1.} implies the statement.  \end{proof}}

\section{The Independence Lemma and How To Use It}\label{sec:indep}

The central result of this section is Proposition \ref{lem3.3} which states that the absence of an $\{a,b\}$-symmetric  polymorphism implies that all polymorphism are canonical on $\{a,b\}$.
The main ingredients of our proof of this proposition are the fact that $\mathbf{A}\in \rra$ has a flexible atom and the following ``Independence Lemma'' (Lemma~\ref{symmetric+formulaslemma}). 


%
%
%

\subsection{The Independence Lemma }

The following lemma transfers the absence of a special partially canonical polymorphism to the existence of certain relations of arity $4$ that are primitively positively definable in $\mathfrak{B}$. 
A lemma of a similar type appeared 
	as Lemma 42 in an article by \cite{RandomMinOps}.

\begin{lemma}[Independence Lemma]\label{symmetric+formulaslemma}
	Let $\mathfrak{B}$ be a homogeneous structure with finite relational signature. Let $a$ and $b$ be 2-orbits of $\Aut(\mathfrak{B })$ such that $a$, $b$, and $(a\cup b)$ are primitively positively definable in $\mathfrak{B }$. 
	Then the following are equivalent:
	\begin{enumerate}
		\item $\mathfrak{B}$ has an $\{a,b\}$-canonical polymorphism $g$ that is $\{a,b\}$-symmetric with $\overline{g}(a,b)=\overline{g}(b,a)=a$.
		\item  
		For every primitive positive formula $\varphi$ such that $\varphi \wedge a(x_1,x_2)   \wedge  b(y_1,y_2)$ and $\varphi \wedge b(x_1,x_2)  \wedge a(y_1,y_2)$ are satisfiable  over $\mathfrak{B}$, the formula $\varphi \wedge a(x_1,x_2)   \wedge  a(y_1,y_2)$ is also satisfiable  over $\mathfrak{B}$.
		\item For every finite $F\subset B^2$ there exists a homomorphism $h_F$ from the substructure of  $\mathfrak{B}^2$ induced by $F$ to $\mathfrak{B}$ that is $\{a,b\}$-canonical with  $\overline{h_F}(a,b)=\overline{h_F}(b,a)=a$.
	\end{enumerate}

\end{lemma}


\begin{proof}
	The implication from (1) to (2) follows directly by applying the symmetric polymorphisms to tuples from the relation defined by $\varphi$.\medskip
	
	For the implication from (2) to (3) let $F$ be a finite subset of $B^2$. Let  $\{e_1,\ldots, e_n\}$ with $n\in \mathbb{N}$ be an enumeration of $F$.
	To construct $h_F$ consider the formula $\varphi_0$ with variables $x_{i,j}$ for $1\leq i,j \leq n$ that is the conjunction of all atomic formulas $R(x_{i_1,j_1},\ldots, x_{i_k,j_k})$ such that $R(e_{i_1}, \ldots, e_{i_k})$ and $R(e_{j_1}, \ldots, e_{j_k})$ hold in $\mathfrak{B}$. Note that this formula states exactly which relations hold on $F$ in $\mathfrak{B}^2$.
	Let $P$ be the set of pairs $((i_1,i_2),(j_1,j_2))$ such that 
	\begin{align*}
		&(a\cup b )(e_{i_1}, e_{i_2}) \\
		\textup{and~}& (a\cup b )(e_{j_1}, e_{j_2}) \\
		\textup{and~}&( a(e_{i_1}, e_{i_2}) \vee  a(e_{j_1}, e_{j_2}) )\\
		\textup{and~}&( b(e_{i_1}, e_{i_2}) \vee  b(e_{j_1}, e_{j_2}) ).
	\end{align*}
	If we show that the formula 
	$$\psi:=\varphi_0 \wedge \bigwedge_{((i_1,i_2),(j_1,j_2)) \in P}  a(x_{i_1,j_1},x_{i_2,j_2})$$
	is satisfiable by an assignment $\alpha$, we get the desired homomorphism by setting $h_F(e_i,e_j):= \alpha(x_{i,j})$. 
	We prove the satisfiability of $\psi$  by induction over the size of subsets $I$ of $P$. 
	For the inductive beginning consider an element $((i_1,i_2),(j_1,j_2)) \in P$. Without loss of generality we have that $a(i_1,i_2)$ holds. Therefore the assignment $\alpha(x_{i,j}):= e_i$ witnesses the satisfiability of the formula $\varphi_0 \wedge  a(x_{i_1,j_1},x_{i_2,j_2})$. 
	For the inductive step let $I\subseteq P$ be of size $m$ and assume that the statement is true for subsets of size $m-1$. Let  $p_1=((u_1,u_2),(v_1,v_2)) $ and $p_2=((u_1',u_2'),(v_1',v_2')) $ be two elements from $I$. 
	We define the following formula
	$$\psi_0:=\varphi_0 \wedge \bigwedge_{((i_1,i_2),(j_1,j_2)) \in I\setminus\{p_1,p_2\}}  a(x_{i_1,j_1},x_{i_2,j_2}).$$
	Then by the inductive assumption the formulas $\psi_0 \wedge a(x_{u_1,v_1},x_{u_2,v_2})$ and 
	$\psi_0 \wedge a(x_{u_1',v_1'},x_{u_2',v_2'})$ are satisfiable. The assumptions on the elements in $P$ give us that also $$\psi_0 \wedge a(x_{u_1,v_1},x_{u_2,v_2}) \wedge b(x_{u_1',v_1'},x_{u_2',v_2'})$$
	and $$\psi_0 \wedge b(x_{u_1,v_1},x_{u_2,v_2}) \wedge a(x_{u_1',v_1'},x_{u_2',v_2'})$$ are satisfiable; since $a\cup b$ is a primitive positive definable relation we are done otherwise.  But then we can apply the assumption of (2) and get that also $\psi_0 \wedge a(x_{u_1,v_1},x_{u_2,v_2}) \wedge a(x_{u_1',v_1'},x_{u_2',v_2'})$ is satisfiable, which proves the inductive step.

	The direction from (3) to (1) is a standard application of K\"onig's tree lemma. For a reference see for example Lemma 42 in the article by \cite{RandomMinOps}.  \end{proof}

\subsection{Absence of $\{a,b\}$-symmetric Polymorphisms}
  We are now able to prove the main result of this section, which  will be a corner stone in the proof of Theorem \ref{theo:result2}.
 Our proof of this proposition makes use of  a $4$-ary relation $ E_{a,b}$ with the following first-order definition:
 \begin{align*}
 	(x_1,x_2,x_3,x_4) \in E_{a,b} :\Leftrightarrow ((a\cup b)(x_1,x_2)\wedge (a\cup b)(x_3,x_4) \wedge a(x_1,x_2) \leftrightarrow a(x_3,x_4)).
 \end{align*}
It is an easy observation that the $\{a,b\}$-canonical polymorphisms of $\frak B$ are precisely those that preserve the relation  $ E_{a,b}$. By Theorem~\ref{polinv} we get that whenever $ E_{a,b}$ is primitively positively definable in $\frak B$ then all  polymorphisms of $\frak B$ preserve $ E_{a,b}$ and are therefore $\{a,b\}$-canonical.  In the following proof we use  the $4$-ary relations that are provided by the second item of  the Independence Lemma~\ref{symmetric+formulaslemma} to provide a primitive positive definition of $ E_{a,b}$.

\begin{proposition}\label{lem3.3}
Let $\mathfrak{B}$ be a normal representation of a finite integral symmetric  relation algebra with a flexible atom $s$. Suppose that $\frak B$ has a binary injective polymorphism. 
Let $a\not \leq\id$ and $b\not \leq \id$  be two atoms such that $\mathfrak{B }$ has no $\{a,b\}$-symmetric polymorphism. Then all polymorphisms are canonical on $\{a,b\}$.
\end{proposition}

\begin{proof}
	
By Corollary \ref{injectivecanonical poly} there exists a canonical binary injective  polymorphism of $\frak B$. Therefore, for every $a'\in A_0$ the edge $\{a',\id\}$ is red and the maximal symmetric polymorphism $t$  that exists by Lemma~\ref{maximal-symmetric defi} is symmetric on all these edges. Observation~\ref{lemma:a id symm means a} implies that $t$ is injective.
Note that $t$ behaves like a projection on $\{a,b\}$ since there exists no $\{a,b\}$-symmetric polymorphism.

Let $\psi$ be the formula defined as follows:$$\psi(x_1,x_2,y_1,y_2):=\overline{\id}(x_1,y_1)\wedge \overline{\id}(x_1,y_2)\wedge \overline{\id}(x_2,y_1)\wedge \overline{\id}(x_2,y_2).$$
We use $\psi$ to formulate and  prove the following claim:
\medskip

\underline{Claim 1:} a) There exists a formula $\varphi_a(x_1,x_2,y_1,y_2)$ such that

\begin{align*}
	&\varphi_a\wedge\psi(x_1,x_2,y_1,y_2)\wedge a(x_1,x_2)   \wedge  b(y_1,y_2) \textup{~ is satisfiable in } \mathfrak{B },\\
	& \varphi_a \wedge\psi(x_1,x_2,y_1,y_2)\wedge b(x_1,x_2)  \wedge a(y_1,y_2) \textup{~ is satisfiable in } \mathfrak{B },\\
&\varphi_a \wedge a(x_1,x_2)   \wedge  a(y_1,y_2) \textup{~ is not satisfiable in } \mathfrak{B }.
\end{align*}%
b) There exists a formula $\varphi_b(x_1,x_2,y_1,y_2)$ that has the same property with  $a$ and $b$ in exchanged roles.

\medskip
\noindent\textit{Proof of Claim 1. }
			There exists a formula $\varphi_a''$ that witnesses the negation of (2) in the Independence  Lemma (Lemma~\ref{symmetric+formulaslemma}) since $\mathfrak{B}$ does not have an $\{a,b\}$-symmetric polymorphism $h$ with $\bar{h}(a,b)=a$. Let $\varphi_b''$ be the formula that witnesses in the same way the non-existence of  an $\{a,b\}$-symmetric polymorphism $h$ of $\mathfrak{B}$ with $\bar{h}(a,b)=b$.
We define for $c\in \{a,b\}$ the formula $\varphi_c'$ as follows:

\begin{equation}
	\varphi'_c(x_1,x_2,y_1,y_2):=  	\varphi''_c(x_1,x_2,y_1,y_2) \wedge (a\cup b)(x_1,x_2)   \wedge  (a\cup b )(y_1,y_2).
\end{equation}

If  $\varphi_a'$ and $\varphi_b'$   witness a) and b) in Claim 1, we are done. So suppose that they do not. Note that if we have  $c\in \{a,b\}$ and $ d\in \{a,b\}\setminus\{c\}$ such that

\begin{align}
	&	\varphi'_c(x_1,x_2,y_1,y_2) \wedge c(x_1,x_2)   \wedge  d(y_1,y_2)\wedge \id(x_2,y_1)  \textup{~ is satisfiable in } \mathfrak{B } \textup{~and} \label{eq1}\\
	&	\varphi'_c(x_1,x_2,y_1,y_2) \wedge d(x_1,x_2)   \wedge  c(y_1,y_2)\wedge \overline{\id}(x_2,y_1) \textup{~ is satisfiable in } \mathfrak{B }, \label{eq2}
\end{align}
 then $\varphi'_c$ would satisfy the statement about $\varphi_c$  in Claim 1, a) or in Claim 1, b). To see this note that we can apply the injective, maximal symmetric polymorphism $t$ that behaves like a projection on $\{a,b\}$ to the tuples $u$ and $v$ that witness (\ref{eq1}) and (\ref{eq2}). The first tuple satisfies $\overline{\id}(x_1,y_1 )$, $\overline{\id}(x_2,y_2 )$ and $\overline{\id}(x_1,y_2 )$ since ${\id}(x_2,y_1 )$ and $c\neq d$. 
 The second tuple satisfies  $\overline{\id}(x_2,y_1 )$. 
 Then the injectivity of $t$ ensures that the  tuples $t(u,v)$ and $t(v,u)$ witness Claim 1, a) or Claim 1, b).
Figure~\ref{fig: poly} illustrates this situation.

\begin{figure}[t]
	
	\centering
	
	\begin{tikzpicture}[scale=1.6]

		\node () at (0,-0.75) [ fill=blue!12, rounded rectangle, rotate=90 ,minimum width=2.9cm,minimum height=0.2cm] {};
		\node () at (0.5,-0.75) [ fill=blue!12, rounded rectangle, rotate=90 ,minimum width=2.9cm,minimum height=0.2cm] {};

		\node () at (2.3,-0.75) [ fill=blue!12, rounded rectangle, rotate=90 ,minimum width=2.9cm,minimum height=0.2cm] {};

		\node (A) at (0,0) [circle,fill,inner sep=1pt, outer sep=2pt] {};
		\node (B) at (0.5,0)[circle,fill,inner sep=1pt, outer sep=2pt] {};

		\node (D) at (0,-0.5)[circle,fill,inner sep=1pt, outer sep=2pt]  {};
		\node (E) at (0.5,-0.5)[circle,fill,inner sep=1pt , outer sep=2pt]{};
		
		\node (G) at (2.3,0)[circle,fill,inner sep=1pt , outer sep=2pt]{};
		\node (H) at (2.3,-0.5)[circle,fill,inner sep=1pt, outer sep=2pt]  {};
			\node () at (0.75,0)[]  {$)$};
			\node () at (-0.25,0)[]  {$($};
			\node () at (0.75,-0.5)[]  {$)$};
			\node () at (-0.25,-0.5)[]  {$($};
		%

		\node () at (1.5,-0.75)[]  {{\LARGE $=$}};

		
		

		\draw[-,thick, >=latex,orange] (A) to (D);
		\draw[-, thick,>=latex, black] (B) to (E);
		
		\draw[-, thick,>=latex, orange] (G) to (H);

		\node (A') at (0,-1) [circle,fill,inner sep=1pt, outer sep=2pt] {};
		\node (B') at (0.5,-1)[circle,fill,inner sep=1pt, outer sep=2pt] {};
		
		\node (D') at (0,-1.5)[circle,fill,inner sep=1pt, outer sep=2pt]  {};
		\node (E') at (0.5,-1.5)[circle,fill,inner sep=1pt , outer sep=2pt]{};
		
		\node (G') at (2.3,-1)[circle,fill,inner sep=1pt , outer sep=2pt]{};
		\node (H') at (2.3,-1.5)[circle,fill,inner sep=1pt, outer sep=2pt]  {};
			\node () at (0.75,-1)[]  {$)$};
			\node () at (-0.25,-1)[]  {$($};
			\node () at (0.75,-1.5)[]  {$)$};
			\node () at (-0.25,-1.5)[]  {$($};

		\draw[densely dotted, thick] (D) to (A');
		
		\draw[dashed] (E) to (B');
		
		\draw[dashed] (A) to[out=-65, in=65] (A');
		\draw[dashed] (D) to[out=-65, in=65] (D');
		\draw[dashed] (A) to[out=255, in=-255] (D');
		
		\draw [-, thick,>=latex] (A') to (D');
		\draw[-, thick,>=latex,orange] (B') to (E');
		
		\draw[-,thick,>=latex] (G') to (H');

		\draw[dashed] (H) to (G');

		\draw[dashed] (G) to[out=-65, in=65] (G');
		\draw[dashed] (H) to[out=-65, in=65] (H');
		\draw[dashed] (G) to[out=255, in=-255] (H');


		\node () at (-0.37,0){{\large $t$}};

	\node () at (-0.37,-0.5){{\large $t$}};
		\node () at (-0.37,-1){{\large $t$}};
			\node () at (-0.37,-1.5){{\large $t$}};
					
		\node () at (0,-1.9)[]  {\large$u$};
		\node () at (0.5,-1.9)[]{\large$v$};
		\node () at (2.3,-1.9)[]{\large$t(u,v)$};

		\node () at (0.23,-2.25)[rotate=325]  {$\in \varphi_c'$};
		\node () at (0.73,-2.25)[rotate=325]{$\in \varphi_c'$};
		\node () at (2.53,-2.25)[rotate=325]{$\in \varphi_c'$};

	\end{tikzpicture}

\caption{Application of the injective operation $t$ on tuples $u$ and $v$. Orange and black edges correspond to atoms $c$ and $d$, the dotted edge to atom $\id$, and the dashed lines denote $\overline{\id}$.}\label{fig: poly}
\end{figure}

By our assumption that Claim 1 is not satisfied by $\varphi_a'$ and $\varphi_b'$  we conclude that for at least one $c\in \{a,b\}$ it holds that for  $ d\in \{a,b\}\setminus\{c\}$  

\begin{align}&\varphi'_c(x_1,x_2,y_1,y_2) \wedge c(x_1,x_2)   \wedge  d(y_1,y_2)\wedge \id(x_2,y_1)  \textup{~ is satisfiable in } \mathfrak{B } \text{~and} \label{eq: 4} \\
&	\varphi'_c(x_1,x_2,y_1,y_2) \wedge d(x_1,x_2)   \wedge  c(y_1,y_2)\wedge \id(x_2,y_1)  \textup{~ is satisfiable in } \mathfrak{B }.  \label{eq: 5}
\end{align}


We distinguish the following different cases. 

%

\begin{enumerate}
	\item $\varphi'_a$ satisfies a)
	in Claim 1 and (\ref{eq: 4}) and (\ref{eq: 5}) hold for $c=b$ and $ d=a$.
	\item (\ref{eq: 4}) and (\ref{eq: 5}) hold for $c=a$ and $ d=b$ and $\varphi'_b$ satisfies b)
	in Claim 1.
	\item (\ref{eq: 4}) and (\ref{eq: 5}) hold for $c=a$ and $ d=b$ as well as for $c=b$ and $ d=a$. 
	
\end{enumerate}

\underline{Case 1:} 	Consider the following formula $\varphi_b$ with
$$\varphi_b(x_1,x_2,y_1,y_2):= \exists z_1,z_2. ( \varphi_b'(x_1,x_2,x_2,z_1) \wedge \varphi_a'(x_2,z_1,z_2,y_1) \wedge \varphi_b'(z_2,y_1,y_1,y_2) ).$$

We claim that $\varphi_b$ satisfies b) in Claim 1. This proves Claim 1, because  $\varphi_a'$ satisfies a) 
 in Claim 1. To see that $\varphi_b$ fulfills the two satisfiability statements in Claim 1, b)  we can first amalgamate the structure $\frak B_1$ induced (as a substructure of $\fB$) by the elements of  a satisfying tuple for $\varphi_b'$ with the structure $\frak B_2$  induced by the elements of a satisfying tuple for $\varphi_a'$ (see Figure~\ref{fig:formula} for an illustration). We amalgamate these two structures over their common substructure $\frak A$ induced by the variables $x_2$ and $z_1$, with the variable names from the definition of $\varphi_b$. In this amalgamation step all missing edges are filled with the flexible atom $s$. In a second amalgamation step we amalgamate the resulting structure with another copy of the structure $\frak B_1$, but now with the common substructure on the variables $z_1$ and $y_1$ (again with refer to the names used in the definition of $\varphi_b$). As before the missing edges are filled with the flexible atom~$s$.
 Figure~\ref{fig:formula} illustrates the situation. 
It follows from the choice of $\varphi_a'$ and $\varphi_b'$ and the definition of $\varphi_b$ that 
$\varphi_b \wedge b(x_1,x_2)   \wedge  b(y_1,y_2)$ is not satisfiable in $\frak B$. \medskip

\begin{figure}[t]
\centering
\begin{tikzpicture}[scale=1]

		\node (A) at (-1,-1)[circle,fill,inner sep=1pt , outer sep=2pt]{};
	\node (C) at (-1,0) [circle,fill,inner sep=1pt, outer sep=2pt] {};
	\node (D) at (0,1)[circle,fill,inner sep=1pt , outer sep=2pt]{};
	\node (C') at (0,2) [circle,fill,inner sep=1pt, outer sep=2pt] {};
	\node (D') at (-1,3)[circle,fill,inner sep=1pt , outer sep=2pt]{};
		\node (F) at (-1,4)[circle,fill,inner sep=1pt , outer sep=2pt]{};
	
	\node () at (-1.7,4){$x_1$};
	\node () at (-1.7,3)[]{$x_2$};


	\node () at (-1.7,-1)[]{$y_2$};
	\node () at (-1.7,0)[]{$y_1$};
	
		\node () at (0.7,1)[]{$z_2$};
	\node () at (0.7,2)[]{$z_1$};

	
		\node (A1) at (-0.8,-1){};
	\node (C1) at (-0.8,-0.1)  {};
	\node (D1) at (0.2,1){};
	\node (C'1) at (0.2,2)  {};
	\node (D'1) at (-0.8,3.1){};
	
	\node (F1) at (-0.8,4){};

	
		\node (A2) at (-1.2,-1){};
	\node (C2) at (-1.2,0)  {};
	\node (D2) at (-0.2,1.1){};
	\node (C'2) at (-0.2,1.9) {};
	\node (D'2) at (-1.2,3){};
	\node (F2) at (-1.2,4){};


		\draw[draw, thick,dashed](A) to (C);
			\draw[draw, thick](D) to (C);
			
			\draw[draw, thick,dashed](D') to (C');
	\draw[draw, thick,black](D') to (F);

		\draw[draw, dotted, out=45, in=-60](A) to (C');
		\draw[draw, dotted, out=110, in=-110](A) to (D');
		\draw[draw, dotted, out=120, in=-120](A) to (F);
	
		\draw[draw, dotted, out=-45, in=60](F) to (D);
	\draw[draw, dotted, out=-110, in=110](F) to (C);



	

	
	

	\draw	[draw=red, 
	fill=red,  fill opacity=0.2] plot [smooth cycle] coordinates { (A1) (C1)(D1)(D2)(C2)(A2) };
	
	\draw	[draw=blue, 
	fill=blue,  fill opacity=0.2] plot [smooth cycle] coordinates {  (C1)(D1)(C'1)(D'1) (D'2)(C'2)(D2)(C2)};
	
	\draw	[draw=red, 
	fill=red,  fill opacity=0.2] plot [smooth cycle] coordinates {  (C'1)(D'1) (F1)(F2)(D'2)(C'2)};

\end{tikzpicture}

\caption{$\varphi_b$ build from $\varphi_b'$  (red), $\varphi_a'$  (blue), the atom $b$ (black), the atom $a$ (dashed) and the flexible atom~$s$ (dotted). The roles of $a$ and $b$ can be changed.}\label{fig:formula}
\end{figure}

\underline{Case 2:} This case is analogous to Case 1. \medskip

\underline{Case 3:} 
Consider the formula $\varphi_a$ with

%
\begin{align*}
	\varphi_a(x_1,x_2,y_1,y_2):= \exists z.( &\varphi_a'(x_1,x_2,x_2,z) \wedge \varphi_b'(x_2,z,z,y_1) \wedge  \varphi_a'(z,y_1,y_1,y_2) ).
\end{align*}
We show that
$\varphi_a\wedge\psi(x_1,x_2,y_1,y_2) \wedge a(x_1,x_2)   \wedge  b(y_1,y_2)$
is satisfiable in $\mathfrak{B }$. Since (\ref{eq: 4}) holds for $c=a$ and $ d=b$ and since $a$ and $b$ are distinct,  there exists $p_1\in A_0\setminus\{\id\}$ such that

\begin{align*}
	\varphi_a'(x_1,x_2,y_1,y_2) \wedge a(x_1,x_2)\wedge b(y_1,y_2) \wedge \id(x_2,y_1)\wedge p_1(x_1,y_2)
\end{align*}
is satisfiable in $\mathfrak{B }$. Similarly, since (\ref{eq: 4}) holds for $c=b$ and $ d=a$, there exists $p_2\in A_0\setminus\{\id\}$ such that

\begin{align*}
	\text{and~~}&\varphi_b'(x_1,x_2,y_1,y_2) \wedge b(x_1,x_2)\wedge a(y_1,y_2) \wedge \id(x_2,y_1)\wedge p_2(x_1,y_2).
\end{align*}

Note that there are $u_1,\ldots, u_5 \in B$ such that the following atomic formulas hold:
\begin{align*}
	a(u_1,u_2),  p_1(u_1,u_3), s(u_1,u_4), s(u_1,u_5),&\\
	b(u_2,u_3), p_2(u_2,u_4),  s(u_2,u_5),& \\
	a(u_3,u_4),  p_1(u_3,u_5),&\\
	 b(u_4,u_5).&
\end{align*} 
If we choose for the existentially quantified variable $z$ in the definition of $\varphi_a$ the element $u_3$ then the tuple $(u_1,u_2,u_4,u_5)$ satisfies the formula $$\varphi_a\wedge\psi(x_1,x_2,y_1,y_2) \wedge a(x_1,x_2)   \wedge  b(y_1,y_2).$$ By an analogous argument also $\varphi_a\wedge\psi(x_1,x_2,y_1,y_2) \wedge b(x_1,x_2)   \wedge  a(y_1,y_2)$ is satisfiable.
It follows again from the choice of $\varphi_a'$ and $\varphi_b'$ and the definition of $\varphi_a$ that 
$\varphi_a \wedge a(x_1,x_2)   \wedge  a(y_1,y_2)$ is not satisfiable in $\frak B$.
By an analogous definition we can find a formula $\varphi_b$ that satisfies b) in Claim~1. Therefore, we are done with Case~3. Altogether this proves Claim~1.

	\medskip

Let $\varphi_a$ and $\varphi_b$ be the two formulas that exist by Claim 1.
	We define the following formulas
		\begin{align*}
&\varphi'_a(x_1,x_2,y_1,y_2):= \varphi_a(x_1,x_2,y_1,y_2) \wedge (a\cup b)(x_1,x_2)  \wedge (a \cup b)(y_1,y_2)\\
&\varphi'_b(x_1,x_2,y_1,y_2):= \varphi_b(x_1,x_2,y_1,y_2) \wedge (a\cup b)(x_1,x_2)  \wedge (a \cup b)(y_1,y_2).
		\end{align*}

	\begin{figure}[t]
		\begin{center}
			\begin{tikzpicture}[scale=0.75]
				\node (C) at (0,0) [circle,fill,inner sep=1pt, outer sep=2pt] {};
				\node (D) at (0,1)[circle,fill,inner sep=1pt , outer sep=2pt]{};
				\node (C') at (0,2.5) [circle,fill,inner sep=1pt, outer sep=2pt] {};
				\node (D') at (0,3.5)[circle,fill,inner sep=1pt , outer sep=2pt]{};

				\node (c) at (-2,0) [circle,fill,inner sep=1pt, outer sep=2pt] {};
				\node (d) at (-2,1)[circle,fill,inner sep=1pt , outer sep=2pt]{};
				\node (c') at (-2,2.5) [circle,fill,inner sep=1pt, outer sep=2pt] {};
				\node (d') at (-2,3.5)[circle,fill,inner sep=1pt , outer sep=2pt]{};
				
				\node (E) at (-1,5)[circle,fill,inner sep=1pt , outer sep=2pt]{};
				\node (F) at (-1,6)[circle,fill,inner sep=1pt , outer sep=2pt]{};

				\node () at (-0.3,5){$x_2$};
				\node () at (-0.3,6)[]{$x_1$};

				\node (B) at (-1,-1.5)[circle,fill,inner sep=1pt , outer sep=2pt]{};
				\node (A) at (-1,-2.5)[circle,fill,inner sep=1pt , outer sep=2pt]{};
				
				\node () at (-0.3,-1.5)[]{$x_3$};
				\node () at (-0.3,-2.5)[]{$x_4$};

					\node () at (0.8,3.5)[]{$y_1$};
						\node () at (0.8,2.5)[]{$y_2$};
							\node () at (0.8,1)[]{$y_3$};
								\node () at (0.8,0)[]{$y_4$};

					\node () at (-2.8,3.5)[]{$z_1$};
				\node () at (-2.8,2.5)[]{$z_2$};
				\node () at (-2.8,1)[]{$z_3$};
				\node () at (-2.8,0)[]{$z_4$};


				\node (C1) at (0.2,0)  {};
				\node (D1) at (0.2,1){};
				\node (C'1) at (0.2,2.5)  {};
				\node (D'1) at (0.2,3.5){};

				\node (c1) at (-1.8,0)  {};
				\node (d1) at (-1.8,1){};
				\node (c'1) at (-1.8,2.5) {};
				\node (d'1) at (-1.8,3.5){};
				
				\node (E1) at (-0.8,5){};
				\node (F1) at (-0.8,6){};

				\node (B1) at (-0.8,-1.5){};
				\node (A1) at (-0.8,-2.5){};


				\node (C2) at (-0.2,0)  {};
				\node (D2) at (-0.2,1){};
				\node (C'2) at (-0.2,2.5) {};
				\node (D'2) at (-0.2,3.5){};

				\node (c2) at (-2.2,0) {};
				\node (d2) at (-2.2,1){};
				\node (c'2) at (-2.2,2.5)  {};
				\node (d'2) at (-2.2,3.5){};
				
				\node (E2) at (-1.2,5){};
				\node (F2) at (-1.2,6){};

				\node (B2) at (-1.2,-1.5){};
				\node (A2) at (-1.2,-2.5){};



				
				\draw[dotted, out=80, in=-80](A) to (E);
				\draw[ dotted](B) to (E);
				\draw[ dotted, out=100, in=260](B) to (F);
				\draw[dotted, out=130, in=230](A) to (F);

%
%
%
%
%
%
%
%

				\draw	[draw=blue, 
				fill=blue,  fill opacity=0.2] plot [smooth cycle] coordinates { (A1) (B1) (c1)(d1)(d2)(c2)(B2)(A2) };
				
				\draw	[draw=red, 
				fill=red,  fill opacity=0.2] plot [smooth cycle] coordinates {  (c1)(d1)(c'1)(d'1) (d'2)(c'2)(d2)(c2)};
				
				\draw	[draw=blue, 
				fill=blue,  fill opacity=0.2] plot [smooth cycle] coordinates {  (c'1)(d'1) (E1)(F1)(F2)(E2)(d'2)(c'2)};

				\draw	[draw=red, 
				fill=red,  fill opacity=0.2] plot [smooth cycle] coordinates { (A1) (B1) (C1)(D1)(D2)(C2)(B2)(A2) };
				
				\draw	[draw=blue, 
				fill=blue,  fill opacity=0.2] plot [smooth cycle] coordinates {  (C1)(D1)(C'1)(D'1) (D'2)(C'2)(D2)(C2)};
				
				\draw	[draw=red, 
				fill=red,  fill opacity=0.2] plot [smooth cycle] coordinates {  (C'1)(D'1) (E1)(F1)(F2)(E2)(D'2)(C'2)};

				%
				%
				%
				%
				%
				%
				%
				%
				%
			\end{tikzpicture}
			
		\end{center}
		\caption{The formula $\delta$ build from $\varphi_a'$  (red), $\varphi_b'$  (blue), and the flexible atom~$s$ (dotted).}\label{fig:formula2}
	\end{figure}

We also define a formula $\delta$ as follows (see also Figure~\ref{fig:formula2}):
	\begin{align*}
	\delta(x_1,x_2,x_3,x_4):= ~&s(x_1,x_3)\wedge s(x_1,x_4) \wedge s(x_2,x_3) \wedge s(x_2,x_4)\\
	&\wedge \exists y_1,y_2,y_3,y_4.(\varphi_a'(x_1,x_2,y_1,y_2)  \wedge \varphi_b'(y_1,y_2,y_3,y_4) \\&\wedge \varphi_a'(y_3,y_4,x_3,x_4))\\
&\wedge \exists z_1,z_2,z_3,z_4.(\varphi_b'(x_1,x_2,z_1,z_2)  \wedge \varphi_a'(z_1,z_2,z_3,z_4) \\&\wedge \varphi_b'(z_3,z_4,x_3,x_4)).
	\end{align*}

Analogously to Case 1, an amalgam of the structures that are induced by tuples that satisfy $\varphi_a'$ and $\varphi_b'$ shows that the formulas $\delta \wedge a(x_1,x_2)  \wedge b(x_3,x_4) $ and $\delta \wedge b(x_1,x_2)  \wedge a(x_3,x_4)$ are satisfiable in $\mathfrak{B }$. Note that this is possible since we ensured in Claim 1 that there exist tuples that additionally satisfy $\psi$. It also holds that a tuple $x$ that satisfies $\delta$ also satisfies 
	\begin{align}\label{eq: a und b3}
( a(x_1,x_2)  \wedge b(x_3,x_4))  \vee   (b(x_1,x_2)  \wedge a(x_3,x_4)).
\end{align}

 Assume that for a tuple $x$ that satisfies $\delta$ it holds that $a(x_1,x_2)  \wedge a(x_3,x_4)$. Then there exist $y_1,y_2,y_3,y_4$ such that $\varphi_a'(x_1,x_2,y_1,y_2) \wedge \varphi_a'(y_3,y_4,x_3,x_4) $ holds. But this is by the definition of $\varphi_a'$ only possible if $b(y_1,y_2)$ and $b(y_3,y_4)$ hold, in contradiction to $\varphi_b'(y_1,y_2,y_3,y_4)$. 
	The same argument works for proving that $\neg ( b(x_1,x_2)  \wedge b(x_3,x_4))$ holds.
	

We complete the proof with a primitive positive definition of $E_{a,b}$. We have the following primitive positive formula $$	\delta'(x_1,x_2,x_3,x_4):=~\exists y_1, y_2.(\delta(x_1,x_2,y_1,y_2) \wedge \delta(y_1,y_2,x_3,x_4)),
$$
and define $E_{a,b}$ by 
$$	(x_1,x_2,x_3,x_4)\in E_{a,b}~ \Leftrightarrow ~ \delta'(x_1,x_2,x_3,x_4).
$$

	For the forward direction of this equivalence, let $x$ be a tuple from $E_{a,b}$ such that $c(x_1,x_2)$ holds  for $c\in \{a,b\}$.  Let  $ d\in \{a,b\}\setminus\{c\}$ and let $y_1$ and $y_2$ be  two elements from $\frak B$ such that $d(y_1,y_2)$  and $s(x_i,y_j)$ for every $i\in \{1,\ldots,4\}$ and every $j\in\{1,2\}$ holds.
%
%
Such elements exists since in the substructure of $\frak B$ that is induced by  $\{x_1,x_2,x_3,x_4,y_1,y_2\}$ all appearing triangles are allowed by the definition of the flexible atom $s$.
	 The elements $y_1$ and $y_2$ witness that $x$ satisfies the formula $\delta'$.
	 
	 For the other direction assume that a tuple $x$ satisfies $\delta'$. Then there exist  $y_1$ and $y_2$ such that $ \delta(x_1,x_2,y_1,y_2) \wedge \delta(y_1,y_2,x_3,x_4)$ is satisfied. Since we observed that the tuples $(x_1,x_2,y_1,y_2)$ and $(y_1,y_2,x_3,x_4)$  both satisfy (\ref{eq: a und b3}) we can assume that $c(x_1,x_2)$ holds  for $c\in \{a,b\}$. It follows also from (\ref{eq: a und b3}) that $d(y_1,y_2)$ holds for $ d\in \{a,b\}\setminus\{c\}$ and then (\ref{eq: a und b3})  implies  that $c(x_3,x_4)$ holds, which proves the backward direction of the stated equivalence.
 
 \end{proof}

\section{Proof of the Result}\label{sec:proof}

{In this section we prove the main results of this article. We first obtain a dichotomy theorem for a class of CSPs  (Theorem~\ref{thm:main 2}). This is used in combination with the observations in Section~\ref{section:symmra with flex atom} to conclude the proof of Theorem~\ref{theo:result2}. 

}

\begin{theorem}\label{thm:main 2}
	Let $\bf A \in \rra$ be finite integral symmetric and with a flexible atom $s$ and let $A_0$ be the set of atoms of $\bf A$. Then either
	\begin{itemize}
		\item there exists an operation 
		$f\colon A_0^6\rightarrow A_0$ 
		that preserves the allowed triples of $\mathbf{A}$, satisfies 
			$$\forall x_1,\ldots,x_6\in A_0.~f(x_1,\ldots,x_6)\in \{x_1,\ldots x_6\}$$ and satisfies the {Siggers identity} 
		$$\forall x,y,z \in A_0.~ f(x,x,y,y,z,z)=f(y,z,x,z,x,y);$$
		in this case, $\csp(\mathfrak{ B}) $ is in P, or 
		\item $\HSPf(\{\Pol(\mathfrak{B }) \}  )$ contains a 2-element algebra where all operations are projections; in this case, $\csp(\mathfrak{ B}) $ is NP-complete. 
	\end{itemize}
\end{theorem}

\begin{proof}
	Let $\mathfrak{B }$ be a normal representation of $\mathbf{A}$ that exists by Proposition~\ref{prop: flex atom implies normal rep}. 
		The finite boundedness of $\mathfrak{B }$ implies  that $\csp(\mathfrak{B} )$ is in NP.
	Let $\mathfrak{O}$ be the atom structure of $\mathbf{A}$. 
We can assume that $\mathfrak{B}$ has a binary injective polymorphism, because otherwise Proposition \ref{binaryinjectionNPhard} would directly imply the second item. The existence of a binary injective polymorphism implies by Corollary \ref{injectivecanonical poly} the existence of a canonical  binary injective polymorphism $g$.

	If the first item of the theorem is satisfied then the operation $f$ is a Siggers polymorphism of   $\mathfrak{O}$ and the statement follows by the Propositions \ref{prop:canon poly=type struc} and \ref{tractabilityB}.
	

Assume therefore that the first item in the theorem does not hold.
By Corollary \ref{trivialsubalgebraO} there exist elements $a,b\in A_0$ such that the subalgebra of $\Pol(\mathfrak{O})$ on $\{a,b\} $ contains only projections. It holds that $\id \not\in \{a,b\}$, since $g$ is a witness that $\id$ can not be in the domain  of a subalgebra that contains only projections. Since all operations from $\Pol(\mathfrak{O})$ are projections on $\{a,b\} $ there exists no canonical polymorphism of $\mathfrak{B}$ that is $\{a,b\} $-symmetric. 
By Proposition $\ref{partiallysymm}$ there exists also no $\{a,b\} $-symmetric polymorphism of $\frak B$.
Since $\frak B$ has a binary injective polymorphism we can apply  Proposition \ref{lem3.3} and get that all polymorphisms of $\mathfrak{B}$ are $\{a,b\} $-canonical.
The last step is to show that all polymorphisms of $\mathfrak{B}$ behave like projections on $\{a,b\} $.

 Assume for contradiction that there exists a ternary, $\{a,b\} $-canonical polymorphism $m$ that behaves on $\{a,b\} $ like a majority or like a minority. By Corollary \ref{partialto global:majomino} there exists a canonical polymorphism that is also a majority or minority on $\{a,b\} $ (here we use again the existence of an injective polymorphism). This contradicts our assumption that $\Pol(\mathfrak{O})$ is trivial on $\{a,b\} $. We get  that every polymorphism of $\mathfrak{B}$ does not behave  on $\{a,b\} $ as an operation from Post's theorem (Theorem~\ref{post}) and therefore must behave as a projection on $\{a,b\} $ by Theorem \ref{post}.   
Thus, $\HSPf(\{\Pol(\mathfrak{B }) \}  )$ contains a $2$-element algebra whose operations are projections and $\csp(\mathfrak{B}) $ is NP-hard, according to Theorem \ref{theo:hardness}.
  \end{proof}

	{ We can prove the main result.
	\begin{proof}[Proof of Theorem~\ref{theo:result2}] Let $\mathbf{A}$ be as in the assumptions of the theorem and let $\bf A'$ be the finite symmetric integral representable relation algebra that exists by Proposition~\ref{prop:flex implies integral}. 
	Suppose that $\mathbf{A}$ satisfies the first condition of Theorem~\ref{theo:result2}. By Item 3) in  Proposition~\ref{prop:flex implies integral}, we get that  then $\mathbf{A}'$ satisfies the first condition in Theorem~\ref{thm:main 2} and therefore $\csp(\frak B)$ for  the normal representation $\frak B$ of $\bf A'$ is in P.  By Section~\ref{subsec: normal rep csp} we know that $\nsp(\bf A')$ and $\csp(\frak B)$ are polynomial-time equivalent. This, together with the Turing reduction from $\nsp(\bf A)$ to $\nsp(\bf A')$ by Item 1) in  Proposition~\ref{prop:flex implies integral} implies that 
	$\nsp(\bf A)$ is in P. This proves the first part of the theorem. 
	
Assume that $\mathbf{A}$ does not satisfy the first condition of Theorem~\ref{theo:result2}. Item 3) in  Proposition~\ref{prop:flex implies integral} again implies that $\mathbf{A}'$ does not satisfy the first condition in Theorem~\ref{thm:main 2} and therefore $\csp(\frak B)$ for  the normal representation $\frak B$ of $\bf A'$ is NP-complete. As before we get that $\nsp(\bf A')$  is NP-complete and the 
many-one  reduction from $\nsp(\bf A)$ to $\nsp(\bf A')$ by Item 2) in  Proposition~\ref{prop:flex implies integral} implies the NP-hardness of $\nsp(\bf A)$. The containment in NP follows by Item 1) in  Proposition~\ref{prop:flex implies integral}. This concludes the proof of Theorem~\ref{theo:result2}. \end{proof}

}
	
	\begin{corollary}
		For a given finite,  symmetric $\bf A \in \rra$ with a flexible atom it is decidable which of the two items in Theorem~\ref{theo:result2} holds.  
		In particular, it is decidable whether $\nsp(\bf A)$ is solvable in polynomial time.
	\end{corollary}
\begin{proof}
	Since $A_0$ is a finite set one can go through all possible operations  $f\colon A_0^6\rightarrow A_0$ that preserve all the allowed triples of $\bf A$ and check whether the Siggers identity is satisfied.
	If  $\text{P}\not=\text{NP}$,  it follows from  Theorem~\ref{theo:result2} that it is decidable whether $\nsp(\bf A)$ is in P. In the case of  $\text{P}=\text{NP}$ this is a trivial task.   
	\end{proof}

The problem of deciding whether certain identities hold in the polymorphism clone of a structure is well known problem in the study of CSPs. The computational complexity is known to be in NP \cite{MetaChenLarose}. The precise complexity of deciding whether a polymorphism clone of an explicitly given structure has a conservative operation that satisfies the Siggers identity is open \cite{MetaChenLarose}.

{

	\section{Connection to Smooth Approximations}\label{sec:approx}
	We discuss the relationship of our results to the techniques developed by \cite{MottetPinskerSmooth}.
	The main invention of  \cite{MottetPinskerSmooth} are \emph{smooth approximations} which are equivalence relations on sets of $n$-tuples. The purpose of  these equivalence relations is to \emph{approximate} the prominent \emph{orbit-equivalence relation}. We have seen the importance of these relations in the present article: the polymorphisms which preserve the orbit-equivalence relation are precisely the canonical ones and they store the information about possible finite-domain algorithms that can be used to solve the infinite-domain CSP (cf. Section~\ref{section:finite type structure}).

	We start by rearranging  the results of Section~\ref{sec:indep} such that we get the following theorem.
	In order to repeat the key steps of our main proof  with a focus on the similarities to  \cite{MottetPinskerSmooth} the assumptions in this theorem are a natural starting point.
	
	\begin{theorem}\label{theo:loop}
		Let $\mathfrak{B}$ be a normal representation of a finite integral symmetric  relation algebra with a flexible atom $s$ and let $\mathfrak{O}$ be the atom structure of $\mathfrak{B}$. Suppose that $\Pol(\frak B)$ contains a binary injective polymorphism and  $\Pol(\mathfrak{O})$ does not have a Siggers operation. Then there exists two atoms  $a\not \leq\id$ and $b\not \leq \id$  such that one of the following holds:
		
		\begin{enumerate}
			\item The orbit-equivalence relation $E_{a,b}$  is primitively positively definable in $\mathfrak{ B}$.
			\item 	For all $x,y\in \{a,b\}$ and every primitive positive formula $\varphi$ such that $\varphi \wedge x(x_1,x_2)   \wedge  y(y_1,y_2)$ and $\varphi \wedge y(x_1,x_2)  \wedge x(y_1,y_2)$ are satisfiable  over $\mathfrak{B}$, the formula $\varphi \wedge x(x_1,x_2)   \wedge  x(y_1,y_2)$ is also satisfiable  over $\mathfrak{B}$.
		\end{enumerate}
		
	\end{theorem}

	This theorem is similar to the ``Loop lemma of approximations'' (Theorem 10 in~\cite{MottetPinskerSmooth}), even though it does not make proper use of the approximation  idea ($E_{a,b}$ approximates itself). The first case in the theorem leads  to the hardness condition 
	that $\HSPf(\{\Pol(\mathfrak{B }) \}  )$ contains a $2$-element algebra whose operations are projections and therefore $\csp(\frak{B})$ is NP-complete.
	
	 However, the second case is ``stronger'' than the second case in Theorem 10 by \cite{MottetPinskerSmooth}. The strength in the statement relies on the special class of problems in our article.
	To see what we mean by this, we continue with the strategy of \cite{MottetPinskerSmooth}.
	As a next step we can restrict the ``Independence Lemma'' (Lemma~\ref{symmetric+formulaslemma})
	as follows:

	\begin{lemma}
		Let $\mathfrak{B}$ be a homogeneous structure with finite relational signature. Let $a$ and $b$ be 2-orbits of $\Aut(\mathfrak{B })$ such that $a$, $b$, and $(a\cup b)$ are primitively positively definable in $\mathfrak{B }$. 
		
		Assume that for every primitive positive formula $\varphi$ such that $\varphi \wedge a(x_1,x_2)   \wedge  b(y_1,y_2)$ and $\varphi \wedge b(x_1,x_2)  \wedge a(y_1,y_2)$ are satisfiable  over $\mathfrak{B}$, the formula $\varphi \wedge a(x_1,x_2)   \wedge  a(y_1,y_2)$ is also satisfiable  over $\mathfrak{B}$. Then 
		$\mathfrak{B}$ has an $\{a,b\}$-canonical polymorphism $g$ that is $\{a,b\}$-symmetric with $\overline{g}(a,b)=\overline{g}(b,a)=a$.
	\end{lemma}
	
	Note that the assumptions in this lemma are precisely what we get from case 2) in Theorem~\ref{theo:loop}. A lemma of similar style can be  also found as Lemma 13 in the article of \cite{MottetPinskerSmooth}. They obtain in this lemma a \emph{weakly commutative} operation. The $\{a,b\}$-symmetric operations from our article are a special case of weakly commutative operations.  The difference between these two properties seems crucial for the next step of  our proof: while $\{a,b\}$-symmetric operations can often be ``lifted'' to canonical $\{a,b\}$-symmetric operations this seems not clear in general for weakly commutative operations. This last step is necessary in order to obtain a contradiction to our assumption that $\Pol(\mathfrak{O})$  does not have a Siggers operation.
	
	In order to apply the results of \cite{MottetPinskerSmooth} directly to obtain the dichotomy results of our article one would have to find a way to canonize the weakly commutative operations in a suitable way. {It seems that this can be done  by a similar proof as for Lemma 56 of  \cite{MottetPinskerSmooth} were the authors used among other things  arguments which are also incorporated in  our proof of the loop lemma (Theorem~\ref{theo:loop}).
	 }

	
}

\section{Conclusion}

We classified the computational complexity of the network satisfaction problem for finite symmetric $\mathbf{A} \in \rra$ with a flexible atom and obtained a P versus NP-complete dichotomy. 
We gave a decidable criterion for $\mathbf{A}$ that is a sufficient condition for the membership of $\nsp(\bf A)$ in P, which is also necessary unless P=NP. 
We want to mention that if we drop the assumptions on $\bf A$ to be symmetric and to have a flexible atom then the statement of Theorem \ref{theo:result2} is false.
An example for this is the Point Algebra; even though the $\nsp$ of this representable relation algebra is in P \cite{PointAlgebra}, the first condition of Theorem~\ref{theo:result2} does not apply.
{However, if we only drop the symmetry assumption we conjecture that Theorem~\ref{theo:result2} still holds; yet it is not clear how to obtain the results of Section~\ref{sec:canon} in this case.} Similarly, if we only drop the flexible atom assumption we conjecture
that the statement also remains true. {For this generalization it would be necessary to obtain Ramsey-type results like Theorem~\ref{theorem:ramsey order} without the assumption of the existence of a flexible atom.}

\subsection*{Acknowledgements} 
The first author has received funding from the European Research Council under the European Community's Seventh Framework Programme (FP7/2007-2013 Grant Agreement no. 681988, CSP-Infinity). The second author is supported by DFG Graduiertenkolleg 1763 (QuantLA).

\bibliography{local}
\bibliographystyle{alpha}

\end{document}